\documentclass{article}
\usepackage[utf8]{inputenc}

\usepackage{fancyhdr}

\usepackage{amsmath}
\usepackage{amsthm}
\usepackage{amssymb}
\usepackage{comment}
\usepackage{soul, color}
\usepackage{appendix}
\usepackage{tasks}
\usepackage{wrapfig}
\usepackage{appendix}
\usepackage{graphicx}
\usepackage[utf8]{inputenc}
 \usepackage[small]{titlesec} 

\usepackage{hyperref}
\usepackage{lipsum}
\usepackage{footnote}
\usepackage{tablefootnote} 
\makesavenoteenv{table}

\newbox{\bigpicturebox}

\usepackage{array}
\usepackage{makecell}

\usepackage{amsmath}

\numberwithin{equation}{section}

\usepackage[ruled,vlined]{algorithm2e}
\usepackage{subcaption,graphicx}

\usepackage{bm}
\usepackage{graphicx}
\graphicspath{{pictures/}}
\DeclareGraphicsExtensions{.pdf,.png,.jpg,.gif}

\usepackage{stackengine} 
\newcommand\oast{\stackMath\mathbin{\stackinset{c}{0ex}{c}{0ex}{\ast}{\bigcirc}}}

\usepackage{makecell}

\usepackage[column=0]{cellspace}
\newcommand{\pd}{\partial}

\newcommand{\N}{\ensuremath{\mathbb{N}}}
\renewcommand{\S}{\mathcal{S}}
\renewcommand{\L}{\mathcal{L}}

\newcommand{\Z}{\ensuremath{\mathbb{Z}}}
\newcommand{\R}{\ensuremath{\mathbb{R}}}

\newcommand{\al}{\alpha} 

\newcommand{\x}{\bm{x}}

\newcommand{\SD}{\Sigma \Delta}

\newcommand{\lam}{\lambda}

\newcommand{\bb}[1]{\mathbf{#1}}

\newcommand{\im}{\mathrm{i}}
\newcommand{\dx}{\mathrm{d}}
 \newcommand{\da}[1]{\partial^{#1}}
\newcommand{\e}{\mathrm{e}}

\newcommand{\ds}{\displaystyle}

\newcommand{\bO}{\mathcal{O}}
\newcommand{\F}{\mathcal{F}}
\newcommand{\B}{\mathcal{B}}

\newcommand{\W}{ \bb W}

\newcommand{\sign}{\mathrm{sign}}

\newcommand{\lbm}{\left|}
\newcommand{\rbm}{\right|}

\newcommand{\lb}{\left(}
\newcommand{\rb}{\right)}

\newcommand{\lbc}{\left\{}
\newcommand{\rbc}{\right\}}

\newcommand{\nofty}[1]{\left\|{#1}\right\|_{\infty}}
\newcommand{\noone}[1]{\left\|{#1}\right\|_{1}}

\usepackage[shortlabels]{enumitem}

\usepackage{relsize}

\usepackage{lipsum} 
\usepackage[OT2,T1]{fontenc}
\DeclareSymbolFont{cyrletters}{OT2}{wncyr}{m}{n}
\DeclareMathSymbol{\Sha}{\mathalpha}{cyrletters}{"58}
 \theoremstyle{plain}
 \newtheorem{thm}{Theorem}[section]
 \newtheorem{lem}{Lemma}[section]

 \newtheorem{prop}{Proposition}[section]

\newtheorem{optprob}{Optimization Problem}[section]

 \theoremstyle{definition}
 \newtheorem{defn}{Definition}[section]
 
 \newtheorem{exmp}{Example}[section]
 \newtheorem{rem}{Remark}[section]

\theoremstyle{remark}

\title{Enhanced Digital Halftoning via Weighted Sigma-Delta Modulation}
\author{ Felix Krahmer, Anna Veselovska\\[6pt]
\small  Technical University of Munich\\ 
\small  Department of  Mathematics
and  Munich Data Science Institute }
\date{ }

\begin{document}
\sloppy

\maketitle

\vspace{-5mm} 

\begin{abstract}


In this paper, we study error diffusion techniques for digital halftoning from the perspective of $1$-bit $\SD$ quantization.
We introduce a method to generate  $\SD$ schemes for two-dimensional signals as a weighted combination of its one-dimensional counterparts and show that various error diffusion schemes proposed in the literature can be represented in this framework via $\SD$ schemes of first order. Under the model of two-dimensional bandlimited signals, which is motivated by a mathematical model of human visual perception,  we derive quantitative error bounds for such weighted $\SD$ schemes. We see these bounds as a step towards a mathematical understanding of the good empirical performance of error diffusion, even though they are formulated in the supremum norm, which is known to not fully capture the visual similarity of images.

Motivated by the correspondence between existing error diffusion algorithms and first-order $\SD$ schemes, we study the performance of the analogous weighted combinations of second-order $\SD$ schemes and show that they exhibit a superior performance in terms of guaranteed error decay for two-dimensional bandlimited signals. In extensive numerical simulations for real world images, we demonstrate that with some modifications to enhance stability this superior performance also translates to the problem of digital halftoning.
More concretely, we find that certain second-order weighted $\SD$ schemes exhibit competitive performance for digital halftoning of  real world images in terms of the {\it Feature Similarity Index} (FSIM), a state-of-the-art measure for image quality assessment. 

\end{abstract}


\section{Introduction}

\subsection{State of The Art}

\begin{figure}[ht!]
\begin{center}
  \subcaptionbox{ 
  \centering
  \label{fig12:a}}{\includegraphics[width=2.1in]{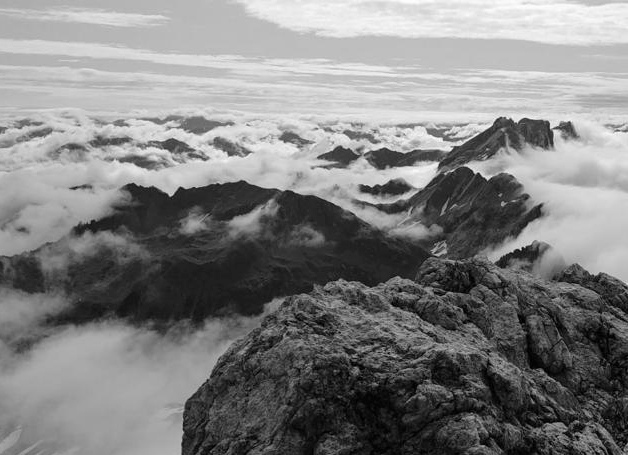}}
 \hspace{10mm}
  \subcaptionbox{  
  \label{fig12:b}}{\includegraphics[width=2.1in]{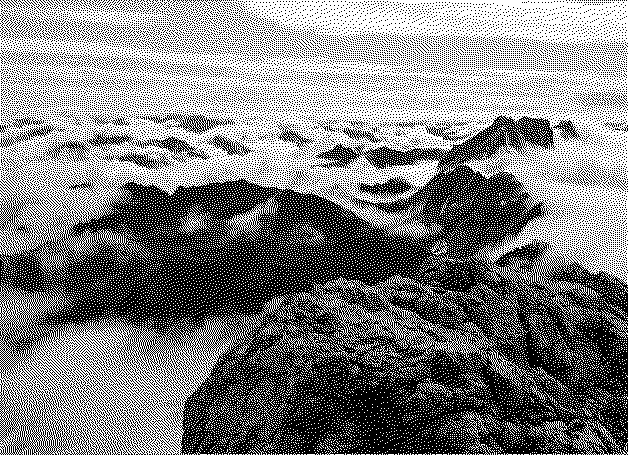}}%
 \caption{Illustration of digital halftoning:
 (a) the original gray-scale image,  (b) the same image represented by black and white pixels using the Floyd–Steinberg algorithm.}
   \label{fig:cat:1}
   
   \end{center}
\end{figure}

Halftoning is an image reproduction technique that simulates continuous-tone imagery through the use of dots. Usually, one differentiates between halftoning or analog halftoning and digital halftoning.  Analog halftoning is a process that simulates  shades of gray or colors by arranging tiny  black, resp. multicolor,  dots of varying size in a regular pattern. The long history of this technique goes back  to  1869 when it was first used in the publishing industry. On the contrary, digital halftoning is a rather modern variant originating in  the 1970s  for use in digital image processing  where the dots 
are of equal size and constrained to the pixel grid. 
An example is given in Figure~\ref{fig:cat:1}: the picture on the left is a gray-scale image while its counterpart on the right is composed of black and white pixels arranged to visually resemble the former as a gray-scale image. 
The key observation that makes it possible is the fact that the human eye acts as a low-pass filter when perceiving visual information from a sufficient distance, blending fine details and recording the overall intensity.
Applications of digital halftoning include not only printing but also sampling problems occurring in rendering \cite{Ostromoukhov2008}, re-lighting \cite{Kollig2003} or object placement
and artistic non-photorealistic image visualization \cite{ Secord2002,Balzer2009}.

In the last decades, methods for digital halftoning developed substantially from basic thresholding and ordered dithering methods to  more sophisticated approaches such as structure-aware halftoning \cite{Pang2008}, which uses optimization techniques to maximize similarity indices, and error diffusion methods, which compute the discrete representation via a recurrence relation. 
Among these state-of-the-art methods, error diffusion techniques are often preferred as they are simpler to implement yet competitive in terms of performance; they will also be the method of choice in this paper. 

There is a significant body of applied literature devoted to the design and analysis of error diffusion techniques.  The most popular error diffusion algorithms are the Floyd-Steinberg algorithm \cite{SteinbergFloyd} and its extensions \cite{Knuth1987, Jarvis1976, Shiau1996}  in which the halftoned image is computed via recurrence relation with fixed coefficients that are not varying spatially or as a function of the gray value.  
More recently, it was discovered in \cite{Ostromoukhov2001} that the Floyd-Steinberg algorithm can give rise to  disturbing patterns at certain gray values. As a solution, the paper proposes a modified  error diffusion algorithm with distribution masks dynamically chosen as a function of the gray value.

In contrast to analog halftoning whose rigorous 
mathematical analysis has been pursued in a number of recent works \cite{Teuber2011, Fornasier2016}, 
the mathematical understanding of digital halftoning in general and error diffusion, in particular, is only in its beginning. 
The goal of this paper is to work towards filling this gap and providing a rigorous mathematical analysis of error diffusion techniques. 

The starting point of our analysis is the connection between error diffusion and $\SD$ quantization that has been hinted at in \cite{Knox1992} and further explored in \cite{Kite1997}. 
$\SD$ quantization  was originally introduced as an analog-to-digital conversion scheme for univariate bandlimited signals \cite{Inose1963, DeFreitas1974}.  
By now, $\SD$ quantization is well-explored in the engineering literature \cite{Aziz1976, Schreier2005}, and initiated by the seminal work of Daubechies and DeVore \cite{Daubechies2003}, its theoretical underpinnings have became an active field in applied mathematics. 
 Some results of relevance to this work will review in Section~\ref{oneD-SD-sec} below.  
 
 Despite the analogies between $\SD$ and halftoning made in \cite{Knox1992, Kite1997}, the error bounds  derived in these works do not directly apply to digital halftoning. The reason is that each pixel of the two-dimensional image needs to be represented by just one bit -- black or white --, while reconstruction guarantees for 1-bit $\SD$ are only available in one dimension.  At the same time, it should be noted that these guarantees are typically formulated in terms of the supremum norm or the mean square error of the underlying bandlimited function, which in our case is the low-pass filtered signal modeling the visual perception. It is not expected that this provides an ideal measure for capturing perceived image similarity. For this reason, we will complement our theoretical analysis in terms of the supremum norm by an in-depth empirical evaluation in terms of state-of-the-art similarity indices.

\subsection{Our Contribution}

In this work, we analyze a class of two-dimensional quantization schemes that arise from weighting and averaging classical 1-bit $\SD$ quantization applied in different directions and show their relevance as error diffusion algorithms for digital halftoning.
On the one hand, we show that a number of error diffusion schemes proposed in the literature can be cast in this framework. On the other hand, we use this approach to design novel algorithms with improved performance. 

To quantify the performance of our algorithms based on weighted $\SD$ schemes and to compare them with other approaches in the literature we follow two different paradigms.

Firstly, we analyze the algorithms' performance as $\SD$ quantization schemes acting on two-dimensional bandlimited signals. This approach measures the performance by comparing the 
the low-pass filtered images as acquired by the human visual system in terms of the supremum norm. While this error metric does not fully align with the visual similarity, it has the advantage of allowing for a rigorous mathematical analysis, giving rise to a quantitative performance measure. This measure can then be used as guidance for the choice of the diffusion coefficients. 

Secondly, we compare the visual quality of the resulting images in terms of the {\it Feature Similarity Index} (FSIM) \cite{Zhang2011}, a state-of-the-art measure for visual similarity of images. We find that error diffusion algorithms constructed via our method from second-order $\SD$ schemes, an approach designed for best performance on bandlimited signals, also outperform a number of popular error diffusion schemes in terms of FSIM. 

A key challenge when implementing and analyzing our method is that $\SD$ schemes are known to provide accurate reconstructions only if they are stable, that is, if the accumulated error remains bounded. Stability, however, is guaranteed for second-order schemes only when the signal amplitude is bounded and indeed, we find that additional stabilizing modifications are required to avoid rare large errors caused by instabilities. Concretely, we observe that a minimal rescaling has hardly any visual effect on the image, while at the same time preventing the rare instabilities due to the underlying second-order $\SD$ schemes.

The paper is organized as follows. We begin with
discussing the mathematical framework in Section~\ref{problem-disc-sec} and fix notation for the
rest of the paper.

Section \ref{oneD-SD-sec} reviews $\SD$ quantization in $1D$ and points out the obstacles in the way of its generalization to the two-dimensional case.  
In Section  \ref{weight-SD-sec}, we introduce and analyze  1st-order weighted $\SD$ quantization schemes; in particular, we establish stability, examine error bounds, and discuss the optimal choice of weights for the supremum norm error metric. 
Weighted $\SD$ quantization schemes of higher order are defined and analyzed in  Section \ref{weight-SD-sec}. 
In Section \ref{numerics-sec}, we confirm the validity of our approach by numerical experiments and explore rescaling as a measure to enhance stability. 
Our conclusions and proposals for future work are discussed in the final section.

\section{Notation and Problem Setting}\label{problem-disc-sec}

The goal of this paper is to mathematically analyze and enhance error diffusion schemes for digital halftoning.  In this section, we discuss the underlying mathematical image models, the digital halftoning problem and the corresponding quantization problem, and different metrics to quantify the representation quality.
At the end of this section, we will also discuss some technical tools necessary for introducing weighted $\SD$ schemes.

\subsection{Image Models, Quantization and Digital Halftoning}

The starting point of our theoretical investigations
is the aforementioned observation that human visual perception involves a smoothing step which can be modeled as a low-pass filter. This property of the human perceptual system together with the fact the scenes observed are not discrete makes the class of bandlimited functions of two variables a suitable  model for visually perceived images.

More precisely, one defines for a bounded region ${\mathcal{R}\subset \R^2}$  the class $\B_{\mathcal{R}}$ of $\mathcal{R}$-bandlimited functions to be the set of real-valued continuous functions in $L^\infty(\R^2)$ 
whose Fourier transforms (as distributions) exist and  
vanish outside of the region $\mathcal{R}$. 
 Here, the Fourier transform of a function $f$ is normalized as 
\begin{equation*}
    \F f( \bm \xi)= \int_{\R^2} f(\bm x)\e^{-2\pi\im \bm \xi \cdot \bm x} \dx \bm x
\end{equation*}
 for each $f\in {L}^1(\R^2)$, and extended to the space of tempered distribution in the usual way. 
 
 Our model is that any visually perceived image can be (at least approximately) represented by some function from the class $\B_{S_{\Omega}}$ for $\Omega>0$ large  enough,  where $S_\Omega$ is the square  $[-\frac{\Omega}{2},\frac{\Omega}{2} ]\!\times\![-\frac{\Omega}{2},\frac{\Omega}{2}]\subset\R^2$  in the frequency domain. For simplicity of presentation we will normalize $\Omega=1$ for remainder of this paper.

To produce digital images, continuous scenes observed by the human eye need to be discretized which mathematically can be understood as a sampling process of the considered model functions.
 Via a  well-known generalization of the Shannon sampling theorem \cite{Petersen1962}, a bandlimited function function $f\!\in\!\B_{S_{1}}$  can be reconstructed  from its samples on the lattice $\!\frac{1}{\lambda}\,\Z^2$ with an oversampling rate $\lambda>1$,  and the sampling formula reads as
\begin{equation}\label{sampling2D}
     f (\x)=\frac{1}{\lambda^2}\sum\limits_{\bm n\in \Z^2} f\Big(\frac{\bm n}{\lambda}\Big)\Phi\Big( \x- \frac{\bm n}{\lambda}\Big), 
\end{equation}
where the kernel $\Phi$ is a Schwartz function with the low-pass property 
\begin{equation}\label{kernel-type}
    \F\Phi(\bm \xi)=\left\{\begin{matrix}
    1,& \bm \xi \in S_{\lambda},\\
    0,&\bm \xi \notin S_{\lambda}.\end{matrix}\right.
\end{equation}



To quantize a signal, we need that it is not only bandlimited but also bounded, which  motivated the definition
\begin{equation}
    \B^\mu:= \lbc f\in \B_{S_1}\colon \; \nofty{f}\le \mu \rbc.
\end{equation}

Given some function $f\in  \B^\mu$, we aim to approximately represent it via $q_{\bm n}$ from some quantization alphabet as 
\begin{equation}
  f_{q}(\x)= \frac{1}{\lambda^2}\sum\limits_{\bm n\in \N^2} q_{\bm n}\Phi\Big( \x- \frac{\bm n}{\lambda}\Big), \quad \quad  \x \in \R^2_{+}. 
\end{equation}
In this paper, we are particularly interested in the case of  $1$-bit quantization, where the alphabet has only two elements. This is due to the fact that in digital halftoning (of a gray-valued image) only the two colors black and white are admissible for each pixel. In the mathematical representation, we renormalize and assume that the elements $ q_{\bm n}$ are chosen from the discrete two-element set~${\mathcal{A}=\lbc-1,1\rbc}$.

If $f_q$ approximates the original function $f$, we call the function $ f_{q}$ {\it a 1-bit representative} of $f$, and the array $q=\{q_{\bm n}\}_{\bm n \in \N^2}$ is referred to as a {\it 1-bit sample sequence}. The {\it main goal} of quantization is to construct a sequence $\lbc q_{\bm n}\rbc_{ \bm n\in \N^2}$ in such that, in a suitable sense,
\begin{equation}
    f_{q}\to f, \quad  \quad  \lambda\to \infty. 
\end{equation}
In our mathematical analysis, we focus on the error metric $\|e\|_{L^{\infty}(\R^2_+)}$, where $e$ is the error signal (or error function)  given by $e(\x):= f(\x)- f_{q}(\x),  \; \x \in \R^2_{+}$. 

In analogy to the one-dimensional case, see \cite{Gunturk2003} for details,  the error signal can be decomposed into two terms
\begin{equation}
    e_{f}(\x):= f(\x)- f_{\lam}(\x),  \quad \quad  e_q(\x):= f_{\lam}(\x)-f_q(\x),
\end{equation}
such that $e=e_{f}+e_{q}$ and the function $f_{\lam}$ is defined as 
\begin{equation}
  f_{\lam}(\x)= \frac{1}{\lam^2}\sum\limits_{\bm n\in \N^2} f\Big(\frac{\bm n}{\lam}\Big)\Phi\Big( \x- \frac{\bm n}{\lam}\Big), \quad \quad  \x \in \R^2_{+},
\end{equation}
It is easy to see that the first term $e_{f}$ does not depend on the quantization approach, but only on $f$ and the oversampling rate $\lambda$, whereas the second term $e_{q}$  depends on both $f$ and the  quantization algorithm. As the kernel $\Phi $ is a Schwartz function, the error  $e_{f}$ will decrease quickly away from zero, which is a direct generalization of an estimate in \cite{Gunturk2003}. 
Motivated by these considerations, our analysis will focus on $e_q$, which we will refer to as {\it quantization error}. 

At the same time, it is well-known that the supremum error norm does not fully capture the perceived visual quality of the image representation. 
As alternatives to error metrics based on such function spaces, a variety of so-called image quality assessment indices have been introduced and demonstrated to better capture visual quality. 

To assess the quality of digital halftoning, such a measure has first been used by Pang et al. \cite{Pang2008}. 
More precisely, the authors consider the Structural Similarity Measure (SSIM) by Wang et al. \cite{Wang2004} and propose to employ {\it structure-aware halftoning} via an iterative optimization method that seeks for a combination of white and black pixels maximizing the SSIM between the original images and its halftone version. 
Such methods have shown very competitive performance. However, because of the discrete optimization step, their main limitation is a fairly long execution time for large images.
Independent of this drawback, the work has made a case for using similarity indices to assess the performance of digital halftoning. 
Inspired by this idea, we will also use image quality assessment indices for measuring the performance of our methods.  
That said,  the drastically different patterns pose a particular challenge to such indices and some of them do not capture the similarities, seeing the reference image and the halftoned image as two entirely different pictures. 
In line with the study of Pang et al.~\cite{Pang2008}, we find that one of the successors of the SSIM, the {\it Feature Similarity Index} for image quality assessment (FSIM) \cite{Zhang2011}, currently one of the most successful and influential full-reference image quality metrics, is particularly successful in capturing the quality of halftoned images, which is why we use this measure in our numerical study. 
The FSIM combines two feature maps derived from the phase congruency measure and the local gradients of the reference and the distorted image to assess local similarities between two images, the reference image, and its distorted counterpart.


\subsection{Directional Differences and Convolutions }\label{Dir-diff-sec}

In this section, we recall some important notions and properties  related to the finite difference operator and the convolution in two-dimensions. The concepts discussed in this section will play a key role in defining weighted $\SD$ schemes.

We begin with introducing the finite-difference operator in the bivariate case. Recall that  the (backward) finite difference $\Delta$ operator maps a sequence $v=\{v_n\}_{n \in \Z}$ to the sequence  $\Delta v= \{(\Delta v)_{ n}\}_{n \in \Z}$ with $(\Delta v)_{ n}=v_n-v_{n-1}$. Consequently, the $r$th order finite difference operator $\Delta^r$ is defined via 
\begin{equation}
    (\Delta^r v)_{ n}=\sum_{j=0}^r (-1)^{j}{\tiny {{r}\choose{j}}}   u_{n-j} 
\end{equation}
When $v$ arises by sampling a smooth  function $f$ with step size $h$, $\Delta^r v$ is known to approximate $h^r\!\cdot\!f^{(r)}$, where $ f^{(r)}$ denotes the $r$th derivative of $f$. Similarly, in two dimensions, $r$th
order finite differences in horizontal and vertical directions approximate the (scaled) partial derivatives. Analogously, finite difference operators can also be defined in arbitrary directions. 

\begin{defn} Let $\L\in \R^2$ be a lattice, for a direction $\bb d= (d_1, d_2)\in \L$ and a two-index sequence $v=\{v_{\bm n}\}_{\bm n\in \L}$, the sequence $\Delta^r_{\bb d}v= \{(\Delta^r_{\bb d}v)_{\bm n}\}_{\bm n\in \Z^2}$ of its {\it directional finite (backward) difference} of order $r$ is defined as 
\begin{equation}
    (\Delta^r_{\bb d}v)_{\bm n}=\sum_{j=0}^r (-1)^{j}{\tiny {{r}\choose{j}}}   v_{\bm n-j\bb d}. 
\end{equation}
\end{defn}

Another interpretation of the action of finite differences is through convolution, which can also be defined for general directions. 

\begin{defn}
 Given a one-dimensional filter $h=\{h_j\}_{j\in \Z}$ supported on the first $L$ elements and a two-index sequence $v=\{v_{\bm n}\}_{\bm n\in \Z^2}$, we denote by $h*_{\bb  d}v$  the {\it convolution} of $v$ and $h$ {\it in the direction} $\bb d$ and define it as
\begin{equation}
    (h*_{\bb d }v)_{\bm n}= \sum\limits_{j=1}^Lh_j v_{\bm n -j \bb d}
\end{equation}
\end{defn}
With this notion, the directional finite difference can be expressed as $(\Delta^r_{\bb d}v)=\Delta^r*_{\bb d } v $ where $\Delta$ denotes the sequence given by  $\Delta_0=1$, ${\Delta_1=-1}$, ${\Delta_k=0}$, for all ${k\in \Z\setminus{\{0,1\}}}$, and $\Delta^r:=\Delta*\ldots*\Delta$.

\section{Background on $\SD$ Quantization }\label{oneD-SD-sec} 

\subsection{ $\SD$ Quantization in $1D$}
The problem of quantizing bandlimited functions on the real line has been studied in a number of works over the last decades. In this section, we briefly review previous works on $\SD$ quantization and introduce  a slightly different perspective on the quantization error analysis, which will form the basis for our generalization to the $2D$-case in the next section.

Consider a  univariate bounded and bandlimited function  $f\in \B^\mu$ and  its sample sequence $y$  as  $y_{n}:=f(\frac{n}{\lam})$ for every $n\in \N$ with the oversampling rate $\lambda>1$. We consider the following generalized form of a $\SD$ quantizer as it has  been introduced and studied  in \cite{Gunturk2003}.

\begin{defn}\label{1d-SD-def} For the sequence $y=\lbc y_{n}\rbc_{n \in \N} $, a {\it $1$-bit $\SD$ quantizer}  takes values of $y$ as input and outputs a sequence $q=\{q_n\}_{n\in \N}$ with  $q_{n}\!\in\!\{-1,1\}$ while constructing the solution to the difference equation
\begin{align}
  v_{n}&= (h*v)_{n}+ y_n-q_n \label{1D-SD-def-1}\\ 
  q_n&= \sign\big((h*v)_{n}+y_{n}\big), \label{1D-SD-def-2}
\end{align}
where the feedback filter  $h\in \ell^1$ satisfies $h_n=0$ if $n\le 0$,  the state variable $v_n$ is set to zero for $n<0$,
and the sign function is given as  
\begin{equation}
    \sign(x):= \lbc \begin{matrix} 1,& x>0, \\
    -1, & x\le 0.\end{matrix} \right.
\end{equation}
\end{defn}

Let  $\delta^a$ denote the Kronecker delta sequence situated at the integer $a$. 
For a positive integer $r$ and a sequence $g\in \ell^1$ with $g_n=0$ for $n<0$, the $\SD$ {\it quantizer} \eqref{1D-SD-def-1}-\eqref{1D-SD-def-2} is called {\it of order $r$} as soon as the filter $h$ satisfies the identity \begin{equation}\label{fil-def}
    \delta^0-h=\Delta^r g. 
\end{equation} 
As shown in \cite{Gunturk2003},  \eqref{1D-SD-def-1} can then be expressed in the more classical form 
\begin{equation}\label{DD-def-of-SD}
    u_n=(\Delta^ru)_n+ y_n-q_n
\end{equation}  via the change of variables $u=g*v$. 

In this paper, we focus on feedback filters that are finitely supported, that is, there exists $L\in \N$ such that $h_n=0$ for all  $n>L$. 
The relation \eqref{fil-def} holds if and only if $ \delta^0-h$ has $r$ vanishing moments, see \cite{Gunturk2003, Krahmerthesis}, that is
\begin{equation}\label{moment-cond}
    \sum\limits_{s=0}^L (\delta^0_{s}-h_s) s^k=0, \quad k=0,\ldots,r-1.
\end{equation}

The first rigorous mathematical error analysis for higher-order quantization schemes (i.e., $r>1$)  described by the recurrence relation \eqref{DD-def-of-SD},  was provided in \cite{Daubechies2003}. This paper is not based on a  quantization rule along the lines of  \eqref{1D-SD-def-1}, but rather constructs the quantizer via a nested sequence of sign operations. For such quantizers the authors derived that reconstruction error of an $r$th-order scheme decays  with the oversampling rate $\lam$ at a rate of $\mathcal{O}({\lam}^{-r})$, and optimizing over the parameter $r$ one can even achieve  a rate of $\mathcal{O}({\lam}^{-c\log\lam})$. By considering schemes of the form  \eqref{1D-SD-def-1}-\eqref{1D-SD-def-2}, G\"unt\"urk  achieved an improved error decay rate of the form $\mathcal{O}(e^{-c\,\lam})$ for some small constant $0<c<1$, that correlates to the maximum admissible  signal amplitude. This type of error decay is optimal: it was shown in \cite{Calderbank2002} the corresponding rate with $c=1$ cannot be achieved; in fact, the maximum feasible $c$ has been shown to decrease to zero when the maximum admissible amplitude increases to one \cite{Krahmer2012}.

To derive error bounds in this paper we will use an approach based on Taylor expansions, which is closely related yet somewhat different from the proof strategies of the aforementioned papers, as we feel that it is better suited to describe our approaches to the two-dimensional scenario.  Namely,  the proof strategy in \cite{Daubechies2003} and \cite{Gunturk2003} is based on summation by parts driven by the finite difference operation in \eqref{DD-def-of-SD}. In higher dimensions, this would require a sequential application of finite difference operators,  while we propose an average.  

To illustrate our alternative approach, we will now rederive some key estimates of one-dimensional $\SD$.

Fix a Schwartz function $\Phi$ satisfying 
\begin{equation}\label{low-pass-cond-ker}
    \F\Phi(\xi)=\left\{\begin{matrix}
    1,&  \xi \in [-\lambda,\lambda],\\
    0,& \xi \notin [-\lambda,\lambda],\end{matrix}\right.
\end{equation} 
which  can be seen as a univariate version of the condition \eqref{kernel-type}.  
Then the quantization  error $e_{q}$ as introduced above, can be expressed as 
\begin{equation*}
    \quad e_{q}(x)= f_{\lam}(x)- f_q(x)= \frac{1}{\lam}\sum\limits_{n=0}^\infty \lb f_{n}-q_{n}\rb\Phi\lb x- \frac{n}{\lam}\rb,  
\end{equation*}
With \eqref{1D-SD-def-1} and reindexing, this yields 
\begin{align}
    e_q(x)=\frac{1}{\lam} \sum\limits_{n=0}^\infty v_n\Big(\Phi\Big(x-\tfrac{n}{\lam}\Big)- \sum\limits_{j=1}^L h_j\Phi\Big(x-\tfrac{n+j}{\lam}\Big) \Big). \label{1d-part-1-of-error}
    \end{align}
Applying the $r$th order Taylor expansion at $a_n= x-\frac{n}{\lam} $, one obtains that 
    \begin{align*}
        \Phi\lb a_n \rb- \sum\limits_{j=1}^L h_j\Phi\Big( a_n -\tfrac{j}{\lam}\Big)& = \Phi\left(a_n\right)- \sum\limits_{j=1}^L h_j\Big(\sum\limits_{p=0}^r \frac{\Phi^{(p)}(a_n)}{p!}\left(\tfrac{-j}{\lam}\right)^p+R_{a_n,r}\left(\tfrac{-j}{\lam}\right)\Big)\\
        &= \sum\limits_{j=1}^L h_j\sum\limits_{p=1}^r (-1)^{p+1} \frac{\Phi^{(p)}(a_n)}{p!}\left(\tfrac{j}{\lam}\right)^p-\sum\limits_{j=1}^L h_j R_{a_n,r}\left(\tfrac{-j}{\lam}\right), \\
      & = (-1)^{r+1}\frac{\Phi^{(r)}(a_n)}{r! \, \lam^r} \sum\limits_{j=1}^L  h_j \cdot j^r-\sum\limits_{j=1}^L h_j R_{a_n,r}\left(\tfrac{-j}{\lam}\right), 
    \end{align*}
  where $R_{a_n,r}(t)=\mathcal{O}(t^{r+1})$ denotes the remainder term, and in the last step we used the vanishing moment conditions for ${\delta^0-h}$.  
Thus \eqref{1d-part-1-of-error} yields
    \begin{equation}
     e_q(x)= \frac{1}{\lam} \sum\limits_{n=0}^\infty v_n\Big((-1)^{r+1}\frac{\Phi^{(r)}(a_n)}{r! \, \lam^r} \sum\limits_{j=1}^L  h_j \cdot j^r-\sum\limits_{j=1}^L h_j R_{a_n,r}\left(\tfrac{-j}{\lam}\right)\Big). \label{1d-part-2-of-error}
        \end{equation}
    The first part of the sum can be estimated as follows
    \begin{align*}
    \frac{1}{\lam} \Big|\sum\limits_{n=0}^\infty v_n\,(-1)^{r+1}\frac{\Phi^{(r)}(a_n)}{r! \, \lam^r} \sum\limits_{j=1}^L  h_j \cdot j^r\Big|  &
    \le\frac{1}{\lam^r} \frac{\big|C_{h}\big|}{r!}\Big| \sum\limits_{n=0}^\infty \frac{v_n}{\lam} \, \Phi^{(r)}(a_n)\Big|\\
   & \le\frac{1}{\lam^r} \frac{\big|C_{h}\big|}{r!} \nofty{v} \|\Phi^{(r)}\|_1
\end{align*}
 where for the $r$th-order filter $h$ the {\it filter constant} is defined by 
 \begin{equation}\label{filter-const}
     C_h:=\sum\limits_{j=1}^L h_s s^r.
 \end{equation} 
The second  part of sum  behaves like $\mathcal{O}\lb\lam^{-(r+1)}\rb$, see Appendix \ref{Taylor-Rem-Est-1st-order} for details. 
Combining these bounds, we obtain  in the following proposition. 
\begin{prop}\label{error-bounds-1d-thm} For a function $f\in \B^\mu$ sampled at rate  $\lam>1$, define the sequence $q\in \{-1,1\}^{\N}$ thought the recurrence \eqref{1D-SD-def-1}-\eqref{1D-SD-def-2}.  Then the  error of the rth order quantization scheme \eqref{1D-SD-def-1}-\eqref{1D-SD-def-2} with a feedback filter $h\in \ell^1(\Z)$ can be characterized as 
\begin{equation}\label{error-bound-1d}
    \nofty{f_{\lam}-f_q}\le  \frac{1}{\lam^r} \|v\|_{\infty} \Bigg(\frac{\big|C_{h}\big|}{r!} \,  \|\Phi^{(r)}\|_1+\mathcal{O}\lb{\lam}^{-1}\rb\Bigg),
\end{equation}
where $\Phi$ is a Schwartz function satisfying the low-pass condition \eqref{low-pass-cond-ker}, 
and  $C_h$ is the filter constant. 
\end{prop}
Note that in the classical case when  the sequence $\delta-h$ coincides with the $r$th order finite difference $\Delta^r$ one has $C_{\Delta^r}=r!$, and the error bound \eqref{error-bound-1d}  corresponds (up to the term  $\mathcal{O}(\lam^{-(r+1)})$) to the error estimates obtained in  \cite{Daubechies2003, Krahmerthesis, Deift2011} by means of repeated integration by parts. 

The error bound estimate \eqref{error-bound-1d} is meaningful only if $\nofty{v}$ is bounded, which motivates the notion of stability. We call the $r$th-order $\SD$ {\it quantization scheme} \eqref{1D-SD-def-1}-\eqref{1D-SD-def-2} {\it stable} if there exists  a number $\eta>0$ such that  
\begin{equation}
    \nofty{v}\le \gamma(r), \quad  \quad  \quad  \text{for all} \quad   \quad  \|y\|_{\infty}\le \eta. 
\end{equation}
A sufficient condition for stability is given by the following proposition, which will also be the criterion of choice in the analysis of this paper.  

\begin{prop}\label{stability-1d} \cite{Gunturk2003}
If the filter $h$ satisfies $\|h\|_{1}\le 2$, then the system \eqref{1D-SD-def-1}-\eqref{1D-SD-def-2} is stable. In particular, for each input sequence $y$ with 
\begin{equation}
    \|h\|_{1}+ \|y\|_{\infty}\le 2
\end{equation}
the state variable $v$ satisfies the bound 
$\|v\|_{\infty}\le 1.$
\end{prop}



\subsection{Towards Stable $2D$ $\SD$ Quantization -- Challenges and First Results}

Generalizing the results on $1$-bit quantization discussed in the previous subsection to two dimensions  is somewhat challenging mainly because stability is harder to achieve. To see this, we note that a natural two-dimensional analogy to  \eqref{DD-def-of-SD} is the recurrence relation  $\ds \Delta_1^{r_1}\Delta_2^{r_2}\,v= f-q$, for integers $r_1,r_2 \ge 0$, and $\Delta_1$ and $\Delta_2$  denoting finite difference operators  acting in the vertical and horizontal direction, respectively. However, for any choice of $r_1, r_2\ge 1$ this leads to filter coefficients with too large of an $l_1$-norm. For example, for $r_1=r_2=1$, we obtain a filter $\delta^0-h$ for 
\begin{equation*}
    h=\begingroup 
\setlength\arraycolsep{2pt}
\begin{pmatrix}
      0 & 1\\
      1& -1 \end{pmatrix} \endgroup
\end{equation*} 
with $\noone{h}=3$. Also introducing an auxiliary sequence $g$ in analogy to the one-dimensional case does not help to overcome this obstacle. Hence, the only chance to achieve the $\noone{h}\le 2 $ is to choose either $r_1\!=\!0$ or $r_2\!=\!0$, which corresponds to performing a one-dimensional $\SD$  scheme either row-by-row or column-by-column. While this approach inherits the recovery guarantees from the one-dimensional case, it is considered by the engineering community as sub-optimal for digital image halftoning as it leads to strong artifacts in the direction orthogonal to the direction of quantization \cite{Knox1992}. 

As demonstrated in \cite{Yilmaz2005} this stability obstacle is specific for $1$-bit quantization and can be overcome by using
2-bit $\SD$ quantization schemes. Recently,  it was shown by Wang and Lyu \cite{Wang2020} that $2$-bit $\SD$ can also be used for efficient image encoding: in contrast to our work, however, these encodings are not proposed as halftoned images, but rather an additional decoding step is required.

 
 In the current work, we introduce an alternative way of generalizing $\SD$  quantization schemes to two dimensions in a stable manner, which in contrast to these works allows for $1$-bit  representations and nevertheless produces high quality halftoned images.



\section{Weighted $\SD$ quantization for Bivariate Signals }\label{weight-SD-sec}

In this section, we present the class of  quantization procedures that are key to this paper, the class of weighted $1$-bit $\SD$ schemes. In contrast to the ideas sketched in the previous section, such schemes combine one-dimensional $\SD$ schemes in different directions in an additive rather than a multiplicative way. As we will show this approach is better compatible with stability, while at the same time it is adapted to the $2D$ signal structure, which is important for the use in digital halftoning.


\subsection{1st-Order Weighted $\SD$ Schemes}

\begin{figure}
\begin{center}
	\includegraphics[width=0.5\textwidth]{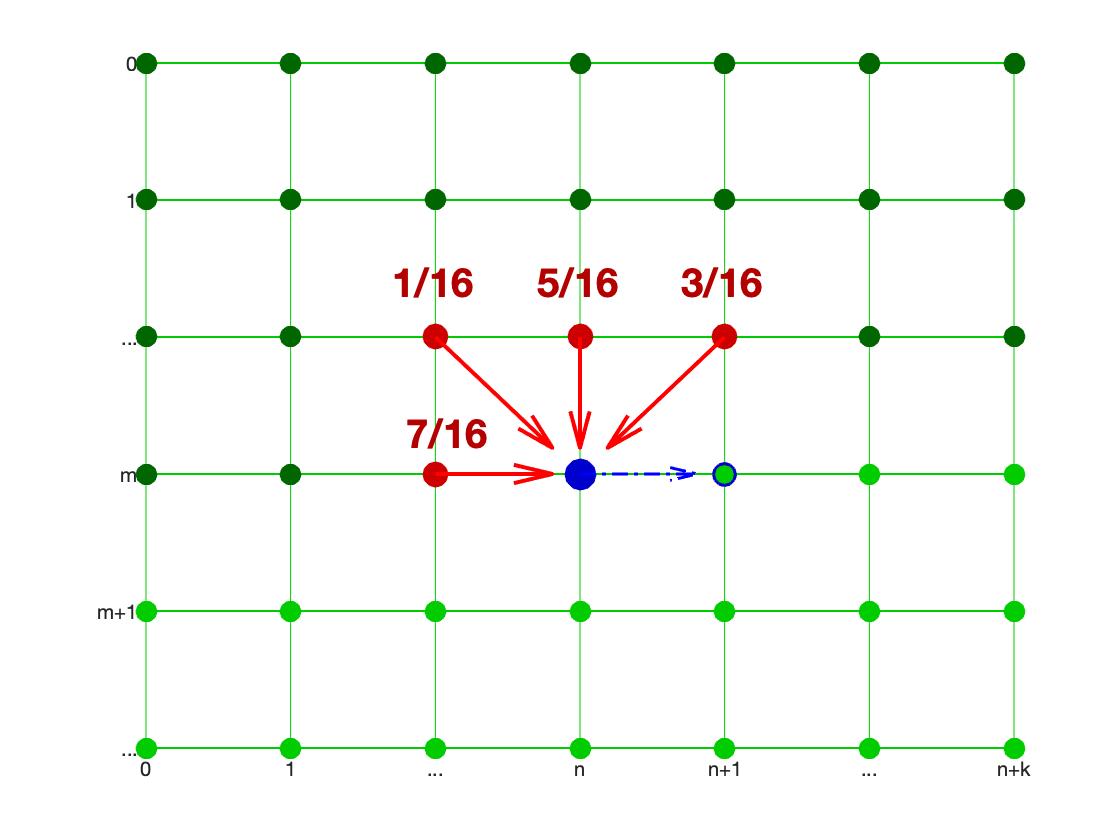} 
	\caption{\small \it The elements of $v$ used (in red) at current quantization step $(m,n)$ (in blue) to define $v_{m,n}$ for  Floyd–Steinberg $\SD$ halftoning scheme. Dark green points denote already half-toned elements and the next step is marked by the green disk with blue bounds. }
		\label{fig:subim1}
		\end{center}
\end{figure}

The main source of inspiration for defining weighted $\SD$ quantization schemes is the celebrated  Floyd–Steinberg halftoning algorithm. 

To describe the algorithm, we consider a gray-scale image as a bivariate sequence  ${y\!=\!\lbc y_{m,n} \rbc_{m, n=0}^{N}}$ of pixel values, which we rescale to the interval $[-1,1]$ for better comparability with  the remainder of this paper. The (analogously rescaled) Floyd–Steinberg algorithm then produces a halftoned image ${q\!\in\!\{-1,1\}^{N\times N}}$ by running the  iterative scheme
\begin{align}\label{SF-scheme-1}
    v_{m,n}&=  \tfrac{5}{16} v_{m-1,n}\!+\!\tfrac{7}{16}v_{m,n-1}\!+\!\tfrac{1}{16}v_{m-1,n-1}\!+\!\tfrac{3}{16}v_{m-1,n+1}\!+\!y_{m,n}-q_{m,n}\\[1pt]
    q_{m,n}&=\sign\Big( \tfrac{5}{16} v_{m-1,n}\!+\!\tfrac{7}{16}v_{m,n-1}\!+\!\tfrac{1}{16}v_{m-1,n-1}\!+\!\tfrac{3}{16}v_{m-1,n+1}\!+\!y_{m,n}\Big). \label{SF-scheme-2} 
\end{align}
see also Figure~\ref{fig:subim1}. 
Comparing this recurrence relation to the one-dimensional $\SD$ quantizer defined in the previous section, we note that can be interpreted as a weighted average of the relation \eqref{DD-def-of-SD}
applied in the four directions ${\bb d_{1,0}\!=\!(1,0)}$, $\bb d_{1,1}\!=\!(1,1),$ $ \bb d_{0,1}\!=\!(0,1),  {\bb d_{-1,1}\!=\!(-1,1)}$. In terms of the directional finite differences introduced in the Section \ref{Dir-diff-sec}, the  Floyd–Steinberg  scheme can hence  be represented  as 
\begin{align*}
    \tfrac{7}{16}(\Delta_{\bb d_{0,1}}v)_{m,n}&\!+\!\tfrac{1}{16}(\Delta_{\bb d_{1,1}}v)_{m,n}\!+\!\tfrac{5}{16}(\Delta_{\bb d_{1,0}}v)_{m,n}\!+\!\tfrac{3}{16}(\Delta_{\bb d_{-1,1}}v)_{m,n} =  y_{m,n}\!-\!q_{m,n}\\ 
q_{m,n}&=\sign\Big( \!v_{m,n}\!-\!\tfrac{7}{16}(\Delta_{\bb d_{0,1}}v)_{m,n}\!-\!\tfrac{1}{16}(\Delta_{\bb d_{1,1}}v)_{m,n}\!-\!\tfrac{5}{16}(\Delta_{\bb d_{1,0}}v)_{m,n}\\
&\quad \quad \quad \quad  \quad  \quad \quad  \quad  \quad \quad  \quad  \quad \!-\! \tfrac{3}{16}(\Delta_{\bb d_{1,-1}}v)_{m,n}\!+\!y_{n,m}\Big)
\end{align*}

\begin{figure}
    \centering
    \includegraphics[width=0.5\textwidth]{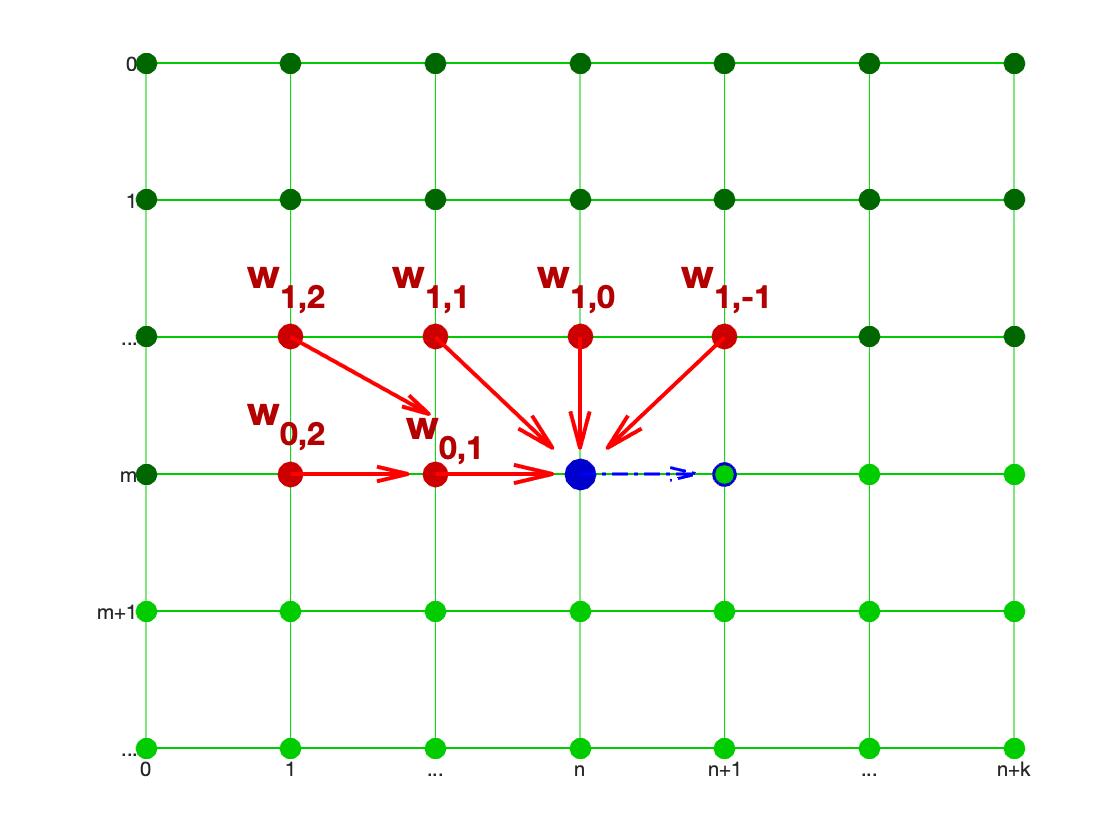}
    \caption{ \it Weighted elements of $v$ (in red) used at current quantization step $\bm n$  to define $v_{\bm n}$ (in blue) for the  weighted  $\SD$ quantization scheme with $\W\in \R^{4\times 2}$. Dark green points denote already quantized elements and the next step is marked by the  green disk with blue bounds.}
    \label{fig:my_label}
\end{figure}

Note that the weights defining the Floyd–Steinberg scheme add up to one and all directions either point to the previous row or a pixel further left in the same row. As we will see these two conditions we be enough to insure stability and also allow for computing the halftoned image via a recurrence relation. 

The first of these two conditions corresponds to combining only $\SD$ schemes in directions $\bb d_{i,j}=(i,j) \in \Z^2$ for which either $i>0$ and $j$ arbitrary or $i=0$ and $j>0$. We restrict our attention to  combinations of such directions and hence assume that  $\bb d_{i,j}$ with $i\in \{0,\dots, p\}$ and $j\in \{-s,\dots,\ell\}$ for some $\ell, s, p \in \N$.   For these $\ell, s, p\in \N$, we average  $\SD$ schemes along different directions with the {\it weight  matrix} $\W \in \R^{(\ell+s+1) \times (p+1)}$ given by 
\begin{equation}\label{weight-m}
    \W= \begingroup 
\setlength\arraycolsep{2pt}
\begin{pmatrix}
     0 &\cdots &0 & w_{0,1}&  \cdots& w_{0,\ell}\\
     \\
    w_{1,-s}& \cdots & w_{1,0}& w_{1,1}&\cdots & w_{1,\ell}  \\
    \vdots& \ddots & \vdots & \vdots &  \ddots & \vdots\\
    w_{p,-s}& \cdots& w_{p,0}& w_{p,1}& \cdots&  w_{p,\ell}
    \end{pmatrix},
\endgroup \quad \text{with} \quad \quad 
     \\
   \sum\limits_{i=0}^{p}\sum\limits_{j=-s}^{\ell}w_{i,j}=1. 
\end{equation}


This gives rise to the following definition. 

\begin{defn} For a given sample sequence $\{y_{\bm n}\}_{\bm n \in \N^2}$,  the  {\it  $1$st-order  weighted $\SD$ quantizer} with weight matrix $\W$ is defined as the iterative scheme 
 \begin{align}
 \label{1st-order-QS-weih_1}
  \sum\limits_{i=0}^{p}\sum\limits_{j=-s}^{\ell}w_{i,j}\left(\Delta_{\bb d_{i,j}} v\right)_{\bm n} &=f_{\bm n}- q_{\bm n}\\[-10pt]
 q_{\bm n}&= \sign\Big( v_{\bm n}- \sum\limits_{i=0}^{p}\sum\limits_{j=-s}^{\ell}w_{i,j}\left(\Delta_{\bb d_{i,j}} v\right)_{\bm n} + y_{\bm n}\Big). \label{1st-order-QS-weih_2}
 \end{align}
\end{defn}

As was mentioned above, the key idea of the 1st-order weighted $\SD$ quantization schemes is to consider a weighted average of $\SD$ quantization schemes of 1st order one applied in different directions. 
As we will see later, when multiple directions are represented by non-zero weights, it helps to smooth out digital halftoning artifacts, on the one hand, and reduce the supremum error for properly chosen weights, on the other hand.  An intriguing fact about the first-order weighted $\SD$ schemes is that they comprise and explain many error diffusion schemes as illustrated in the examples below.

\begin{exmp}\label{exp-halftone} In the following list examples we denote  the element $w_{0,0}$ on the weight matrix by zero in bold, to indicate how many negative directions are included. 

\vspace{-1mm}

\begin{enumerate}[(a), itemsep=0.5pt,parsep=0pt]

\item Applying the one-dimensional  1st-order $\SD$ scheme to the bivariate samples {\it row-by-row} corresponds  the weight matrix

\vspace{-1mm}

\begin{equation*}
    \W_{RbR}= \left(\begin{matrix}
      \bm{ 0} & 1\\
     0& 0  \end{matrix}\right).
\end{equation*}

\vspace{-1mm}


\item {\it Simple averaging} over two perpendicular directions can be represented by  the weight matrix

\vspace{-1mm}

\begin{equation*}
    \W_{1/2}= \left(\begin{matrix}
      \bm{ 0} & \frac{1}{2}\\
      \frac{1}{2}& 0  \end{matrix}\right). 
\end{equation*}


\noindent  We will see that despite its simple structure this scheme exhibits remarkably good performance for digital halftoning. 


\item The {\it  Floyd–Steinberg scheme} has the  the weight matrix 
\vspace{-1mm}

\begin{equation*}
    \W_{F\text{-}S}= \left(\begin{matrix}
     0 & \bm{ 0} & \frac{7}{16}\\[1pt]
     \frac{3}{16} & \frac{5}{16}& \frac{1}{16}    \end{matrix}\right).
\end{equation*}

\item One of  the {\it Shiau-Fan schemes} \cite{Shiau1996}, which were introduced as improvements of the Floyd–Steinberg algorithm  corresponds to the weight matrix 

\vspace{-1mm}

\begin{equation*}
     \W_{Sh\text{-}Fan}= \left(\begin{matrix}
 0&  0 & 0 & \bm{ 0} & \frac{8}{16}\\[1pt]
 \frac{1}{16} &\frac{1}{16}     & \frac{2}{16} & \frac{4}{16}& 0   \end{matrix}\right).
 \end{equation*}


\item { \it The 12-element Jarvis-Judice-Ninke scheme} well-known as an edge enhancement technique \cite{Jarvis1976}  can be represented by  the weight matrix 

\vspace{-3mm}

\begin{equation*}
    \W_{JJN}\!=\! \left(\begin{matrix}
    0& 0 & \bm{ 0} & \frac{7}{48}& \frac{5}{48} \\[2pt]
     \frac{3}{48}& \frac{5}{48}& \frac{7}{48}& \frac{5}{48}& \frac{3}{48}\\[2pt]
     \frac{1}{48}& \frac{3}{48}& \frac{5}{48}& \frac{3}{48}& \frac{1}{48} 
     \end{matrix}\right).
\end{equation*}
\end{enumerate}
\end{exmp}

The following theorem provides an estimate  for the quantization error of 1st-order weighted $\SD$ schemes in terms of the oversampling rate.

\begin{thm}\label{thm-2d-error-bounds-1st-WSD}
Consider a bandlimited function $f\in  \B^\mu$ sampled on the lattice $ \frac{1}{\lam}\,\N^2$ with oversampling rate $\lam>1$. Then  the 1-bit sequence $q\!\in\!\{-1,1\}^{\N^2}$ constructed by the  1st-order weighted $\SD$ quantization scheme  \eqref{1st-order-QS-weih_1}--\eqref{1st-order-QS-weih_2} defines a quantized representative $f_{ q}$ such that 
\begin{equation}\label{error-bound-1st-order-weighted}
    \nofty{f_{\lam} - f_q} \le  \frac{1}{\lam} \nofty{v} \, \Big(   C_{\W} \cdot C \cdot \|\nabla \Phi\|_{1,2}+ \mathcal{O}\big(\lam^{-1}\big)\Big)
\end{equation}
where $C>0$ is a constant independent of $\W$, $\Phi$ is a Schwartz function of the  low-pass type \eqref{kernel-type},  and the absolute constant $C_{\W}$ is determined by 
\begin{equation*}
   \lb C_{\W}\rb^2: = \Big( \sum\limits_{i=1}^p\sum\limits_{j=-s}^{\ell}   i w_{i j}  \Big) ^2  +   \Big( \sum\limits_{i=0}^p\sum\limits_{j=-s}^{\ell}   j  w_{i j}  \Big)^2.
\end{equation*}
\end{thm}

For the proof of Theorem~\ref{thm-2d-error-bounds-1st-WSD} see Section~\ref{sec-error-est-WSD}.

In the estimate \eqref{error-bound-1st-order-weighted}, the leading  error term depends on the {\it weight constant}  $C_{\W}$, and consequently weighted  $\SD$ schemes corresponding  to  $\W$ with small $C_{\W}$ are expected to yield small quantization error. Indeed, the weight constants substantially differ for different weight matrices.  For instance,  the Shiau-Fan approach  provides $C_{\W_{S\text{-}Fan}}\!\approx\!0.5$, while the row-by-row $\SD$ scheme in Example~\ref{exp-halftone} results in $C_{\W_{RbR}}\!=\!1$, see Table~\ref{tab1} for details, and we expect a corresponding gain in the reconstruction error. As our numerical simulations show, these improved constants actually translate to improved reconstruction accuracy. This motivated us to  minimize the weight via the following optimization problem.

\begin{optprob}\label{opt-prob-1st-or} $\quad$
 For fixed $ \; \ell,\; s,\; p\in \N \;$ find 
 \begin{eqnarray*}
  \min\limits_{\W\in \R_{>0}^{(\ell+s+1)\times (p+1)}} \quad   \Big(   \sum\limits_{i=0}^{p}\sum\limits_{j=-s}^{\ell}j w_{i,j} \Big)^2+  \Big(  \sum\limits_{i=1}^{p}\sum\limits_{j=-s}^{\ell} i w_{i,j} \Big)^2 \; \\[-2pt]
    \quad\mbox{ subject to }\quad {  \sum\limits_{i=0}^{p}\sum\limits_{j=-s}^{\ell}w_{i,j}=1},\\
    w_{i,j}\ge 0, \\
    w_{0,-j}=0, \quad  j=0,\dots,s. 
 \end{eqnarray*}
 \end{optprob}
 
\begin{figure}
\begin{center}
	\includegraphics[width=0.5\textwidth]{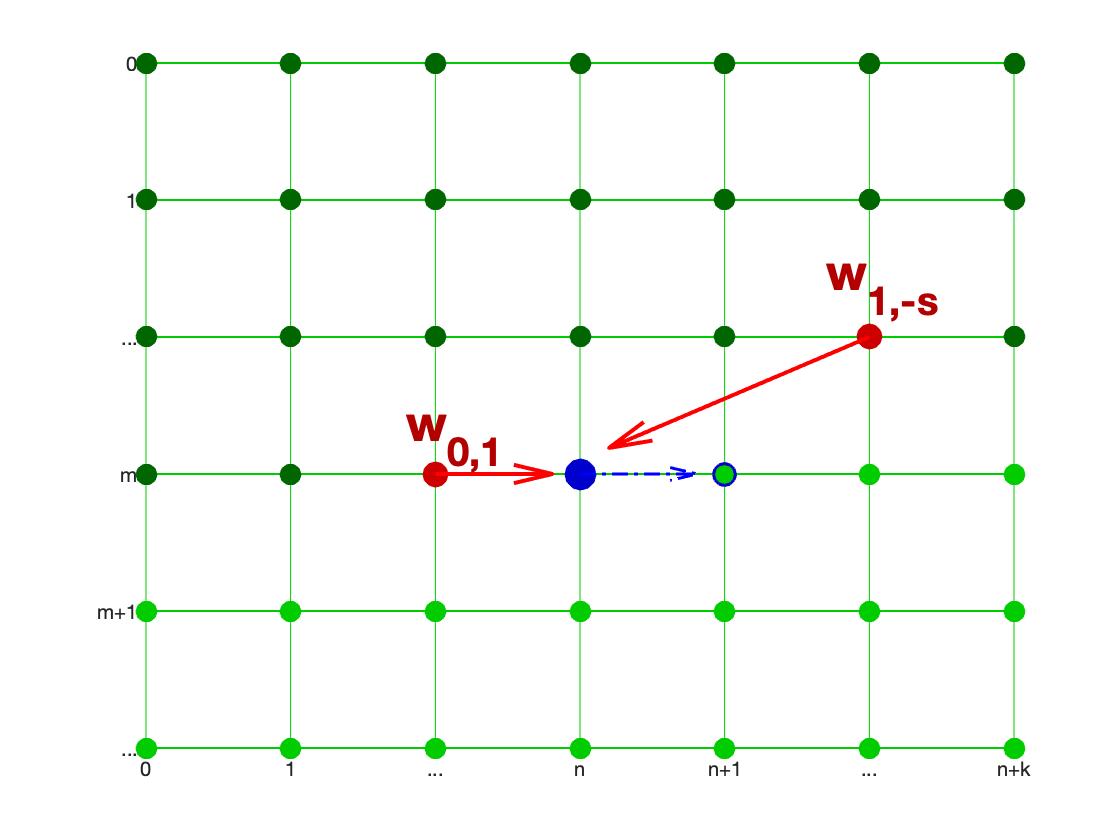} 
	\caption{\small \it Choice of optimal weights for $ \W \in \R^{4 \times 2}$. }
		\label{fig:4}
		\end{center}
\end{figure}

 The solution to this problem  can be explicitly computed, namely its minimal value 
\begin{equation*}
C_{\W_\mathrm{min}}=\frac{1}{\sqrt{1+(s+1)^2}}
\end{equation*}
is attained when the non-zero weights are taken as
\begin{equation*}
  w^{\mathrm{min}}_{1,-s}=\frac{s+1}{1+(s+1)^2} \quad \quad  \text{and} \quad \quad
  w^{\mathrm{min}}_{0,1}=1-\frac{s+1}{1+(s+1)^2},
\end{equation*} 
and all other weight coefficients are set to zero, see Figure~\ref{fig:4} for illustration.


While this very sparse design may seem counter-intuitive, we find in our numerical experiments that for the scenario it is meant to optimize, namely, for two-dimensional bandlimited functions, it indeed outperforms other approaches in terms of the supremum norm. That said, 
this improvement is not reflected in an improved visual quality of the halftoned image as quantified by the similarity index FSIM. In that sense, this example illustrates the discrepancy between sampling theory and halftoning practice; the supremum norm is not an ideal measure for visual quality.  Even though the optimal choices of the weight matrices for the supremum norm and the visual halftoning quality do not agree, however, we find that a number of weighted 1st-order $\SD$ schemes with near-optimal weight constants exhibit very good performance also for digital half-toning. We see this as evidence that the sampling theory perspective yields an error analysis that at least approximately captures the performance of error diffusion, and can hence provide a general idea of why these methods work for digital halftoning. These insights are then of crucial importance for designing weighted $\SD$ schemes of higher order.

\subsection{Weighted $\SD$ Schemes of Higher Order}\label{higher-order-SD-section}

 As discussed above the error diffusion schemes such as  Floyd–Steinberg,  Shiau-Fan, or Jarvis-Judice-Ninke,  can be interpreted as weighted 1st-order quantization schemes.  
 
 Given the superior performance of higher-order $\SD$ scheme in one dimension, it is a natural question whether also weighted $\SD$ schemes of higher order can be used for digital halftoning. To the best of our knowledge,  however, no error diffusion schemes studied in the literature can be interpreted as a weighted higher-order $\SD$ scheme. 
 
 An explanation why none of the ad-hoc error diffusion schemes admit such an interpretation may be the aforementioned observation that to guarantee stability for higher order schemes one typically needs an amplitude less than one, which can in general not be assumed for gray-scale images. 
 
 In this paper, we nevertheless propose weighted  $\SD$  schemes of second order for digital halftoning. The reason is that by choosing appropriate filters one obtains very mild amplitude constraints, which can be addressed by a minimal rescaling with hardly any visual effect.   Thus, in addition to a weight matrix $\W \in \R^{(\ell+s+1) \times (p+1)}$ as given in  \ref{weight-m}, we also need to carefully choose  the feedback filters $h$ of the underlying higher-order $\SD$ schemes in $1D$ (cf.~Definition~\ref{1d-SD-def}).

\begin{figure}
    \centering
    \includegraphics[width=0.7\textwidth]{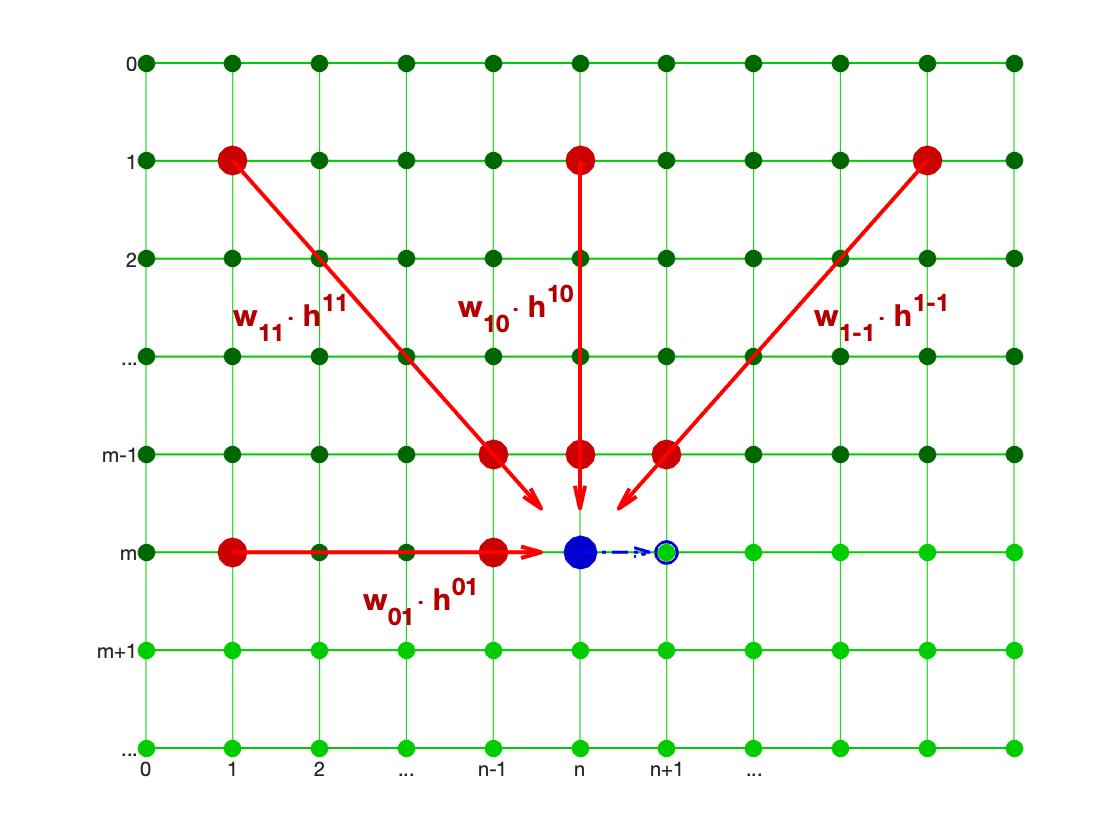}
    \caption{ \it  Visualization of a weighted $\SD$ scheme of second order with weight matrix $\W\in \R^{3\times 2}$ and a  ($2$-sparse) filter $h\in \R^4$. The red dots indicate elements of $v$ used in the current quantization step to define the $v_{\bm n}$  indicated by a blue dot. Dark green dots denote already quantized elements, light green dots those not yet quantized. The subsequent element to be quantized is indicated by a blue circle.}
    \label{fig:my_label-2}
\end{figure}


 We obtain the following definition.

\begin{defn}\label{high-order-sch-def} For a given sample sequence $y=\{y_{\bm n}\}_{\bm n \in \N^2}$,  the  {\it  weighted $\SD$-quantizer} with the weight matrix $\W$ and is given by 
\begin{align}\label{rth-order-sim-QR-filter_1}
  v_{\bm n}-\sum\limits_{i \, j}w_{i,j}\left(h^{i,j}*_{\bb d_{i,j}} v\right)_{\bm n} &=y_{\bm n}- q_{\bm n}\\
 q_{\bm n}&= \sign\Big(   \sum\limits_{i \, j} w_{i,j}\left(h^{i,j}*_{\bb d_{i,j}}v\right)_{\bm n} + y_{\bm n}\Big), \label{rth-order-sim-QR-filter_2}
 \end{align}
 where the single-index filters $h^{i,j}=\lbc h^{i,j}_n\rbc_{n\in \Z}\in \ell^1(\Z)$. If all filters $h^{i,j}$ fulfill the condition $$\delta^0-h^{i,j}= \Delta^r g^{i,j}$$
 for sequences $g^{i,j}\in \ell_1(\Z)$ with  $g^{i,j}_n=0$ for $n<0$,  we speak of a quantizer of  order $r$. 
 \end{defn}
 
 The linear dependence on $v$ in the equations 
 \eqref{rth-order-sim-QR-filter_1}-\eqref{rth-order-sim-QR-filter_2} can be described by a single {\it extended weight matrix} incorporating the effects of both $h$ and $\W$.  We illustrate this in the following example.

\begin{exmp} In all of the following examples,  we consider  the filters ${h^{2}:=[0, \; 3/2, \; 0, \; -1/2]}$ and ${h^{3}=[0, \; 4/3, \; 0, \; 0, \; -1/3]}$ which are particular examples of a larger class of sparse second-order filters defined in \cite{Gunturk2003, Krahmerthesis,Deift2011} that will also be employed in our numerical experiments in Section~\ref{num-halftone-sec}.  

\begin{enumerate}[(a)]
\item  {\it Simple average} schemes of second order use the weight matrix
\begin{equation*}
    \W_{1/2}= \left(\begin{matrix}
      0 & \frac{1}{2}\\
    
      \frac{1}{2}& 0  \end{matrix}\right), 
\end{equation*}

combined with second-order filters $h^{i,j}$. For $h^{0,1}\!=\!h^2$  and $h^{1,0}\!=\!h^3$ or for $h^{0,1}\!=\!h^{1,0}\!=\!h^3$, respectively, we obtain the extended weight matrices 
\begin{equation*}
\W^{2nd}_{A23}= \left(\begin{matrix}
     \bm 1 & -\frac{3}{4} & 0 &  \frac{1}{4}\\
     -\frac{4}{6}&0& 0&0 \\
     0 & 0& 0&0  \\[-2pt]
     0 & 0& 0&0 \\[-1pt]
\frac{1}{6} &0& 0&0 \\
      \end{matrix}\right),  
      \quad \quad 
    \W^{2nd}_{A33}=  \left(\begin{matrix}
     \bm 1 & -\frac{4}{6} & 0 & 0& \frac{1}{6}\\
     -\frac{4}{6}&0& 0&0 &0\\
     0 & 0& 0&0 &0 \\[-2pt]
     0 & 0& 0&0 &0\\[-1pt]
\frac{1}{6} &0& 0&0 &0\\
      \end{matrix}\right), 
\end{equation*}

and as it is easy to see, the matrix $\W_{2nd\text{-}A22}$ describes the $\SD$ scheme with all filters equal to $h^2$, $\W_{2nd\text{-}A33}$ corresponds to the case with all $h^3$-filters, and $\W_{2nd\text{-}A23}$ emerges from the combination of $h^2$ and $h^3$ in the appropriate directions.  

\item The {\it  Floyd–Steinberg schemes of second order} have the weight matrix 
\begin{equation*}
    \W_{F\text{-}S}= \left(\begin{matrix}
     0 & \bm{ 0} & \frac{7}{16}\\[2pt]
     \frac{3}{16} & \frac{5}{16}& \frac{1}{16}    \end{matrix}\right), 
\end{equation*}
so when all four feedback filters $h^{i,j}$ are chosen to be  $h^3$  we obtain the extended weight matrix
\begin{equation*}
    \W^{2nd}_{F\text{-}S\text{-}33}=  \left(\begin{matrix}
    0 & 0&  0 & 0&   \bm 1 & -\frac{28}{48} & 0 & 0& \frac{7}{48}\\[2pt]
   0 & 0&  0 & -\frac{12}{48}&    -\frac{20}{48}&-\frac{4}{48} & 0&0 &0\\
   0 & 0&  0 & 0&    0 & 0& 0&0 &0 \\[-2pt]
   0 & 0&  0 & 0&    0 & 0& 0&0 &0\\[-2pt]
 \frac{3}{48} & 0&  0 & 0& \frac{5}{48} &0& 0&0 &\frac{1}{48} \\
      \end{matrix}\right).
\end{equation*}
In the similar fashion, one can also design the second-order Shiau-Fan and Jarvis-Judice-Ninke schemes. 
\end{enumerate}
\end{exmp}

The next result establishes error bounds for  weighted $\SD$ scheme of higher order; its proof is discussed in Section~\ref{sec-error-est-WSD}.

\begin{thm}\label{thm-2d-error-high-order-WSD} 
Consider a bandlimited function $f\in  \B^\mu$ sampled on the lattice $ \frac{1}{\lam}\,\N^2$ with the oversampling rate $\lam>1$. If a weighted $\SD$ scheme \eqref{rth-order-sim-QR-filter_1}-\eqref{rth-order-sim-QR-filter_2} used for construction of $f$'s 1-bit samples $q\!\in\!\{-1,1\}^{\N^2}$ is of order $r$, then the corresponding quantized representative $f_{ q}$ satisfies 
\begin{equation}\label{error-bound-high-order}
    \nofty{f_{\lam} - f_q} \le  \frac{1}{\lam^r} \,    \nofty{v} \Big( C_{\W}\cdot C \cdot \|\nabla^r \Phi\|_{1,2}  + \mathcal{O}(\lam^{-1})\Big),
\end{equation}
where $C>0$ is a constant independent of $\W$, $\Phi$ is a Schwartz function of the low-pass type \eqref{kernel-type},  the weight constant of order $r$ is defined as 
\begin{equation}\label{const-of-weigths-rth-order}
   (C_{\W})^2:=\sum_{m=0}^r\Big( \sum\limits_{i=0}^{p}\sum\limits_{j=-s}^{\ell} w_{i,j} \cdot C_{h^{i,j}} \cdot i^{r-m}j^m\Big)^2,
\end{equation}
with filter constants $C_{h^{i,j}}$ as in \eqref{filter-const}, and 
   $ \left\|  \nabla^r \Phi\right\|_{1,2}:=  \sqrt{\sum_{|\alpha|=r}\tfrac{1}{\alpha!}\noone{\pd^{\alpha}\Phi}^2}.$ 
   
\end{thm}

Similarly to the 1st-order schemes, one can consider an optimization problem to minimize the weight constant $C_{\W}$. 

 \begin{optprob} For fixed $ \; \ell,\; s,\; p\in \N \;$ find 
   \begin{eqnarray*}
  \min\limits_{\W\in \R^{(\ell+s+1)\times (p+1)}} \quad  \sum\limits_{m=0}^{r} \Big(   \sum\limits_{i=0}^{p}\sum\limits_{j=-s}^{\ell} C_{h^{i,j}} w_{i,j} \,i^{r-m}j^m \Big)^2 \; \\[-3pt]
    \quad\mbox{ subject to }\quad {  \sum\limits_{i=0}^{p}\sum\limits_{j=-s}^{\ell}w_{i,j}=1},\\[-3pt]
    w_{i,j}\ge 0, \\[-3pt]
    w_{0,-j}=0, \quad  j=0,\dots,s.
 \end{eqnarray*}
   \end{optprob}
 Here the optimal weights depend on  the filter constants, and can be numerically computed 
 for each weighted $\SD$ scheme individually. For instance, for the weighted  $\SD$ of 2nd order with all filters $h^{i,j}$ set to $h=[0, \; 4/3, \; 0, \; 0, \; -1/3]$,  the weight constant $C_{\W}$ reaches it minimal value of $\frac{5}{3 \sqrt{2}}$  for  $w_{0,1}\!=\!w_{1,0}\!=\!1/2$ and all other weights equal to zero. At the same time,  for the row-by-row scheme with the same filter, one obtains  $C_{ \W}\!=\!\frac{5}{3}$, see Table~\ref{tab2} for more details.

\subsection{Stability and Error Estimation for Weighted $\SD$ Schemes }\label{sec-error-est-WSD}

As in the one-dimensional case, also the error bounds for weighted $\SD$ schemes in two dimensions are meaningful only if the state variable is bounded. Hence, a key step towards making these results applicable is to establish stability. 



The following proposition generalizes   the rigorous stability analysis in  \cite{Gunturk2003}  to the  two dimensional case providing a sufficient condition for the stability of weighted  1-bit $\SD$  quantization schemes.

\begin{prop}\label{Stability-prop}
For any filter $h^{i,j}$ satisfying  $\|h^{i,j}\|_{1}< 2$ for all $i=1,\dots, p$ and $j=-s, \dots, \ell$,  and initialization $\lbm v_{\bm n} \rbm \le 1$ for $\bm n \in \N^2$, the two-index sequence $v=\{v_{\bm n}\}_{\bm n \in \N^2}$ defined by the recursion  \eqref{rth-order-sim-QR-filter_1}-\eqref{rth-order-sim-QR-filter_2} is uniformly bounded. Namely, for each sample sequence $y$ with 
\begin{equation}
    \sum\limits_{i \, j} w_{i,j}  \noone{h^{i,j}}  + \|y\|_{\infty}\le 2
\end{equation}
the state variable $v$ satisfies the bound 
$\|v\|_{\infty}\le 1.$
\end{prop}

\begin{proof}
We will argue by induction.  Suppose that starting with zero the weighted $\SD$ scheme \eqref{rth-order-sim-QR-filter_1}-\eqref{rth-order-sim-QR-filter_2} is applied $\ell-1$
times. 
Assume $\lbm v_{\bm n_{k}} \rbm \le 1$ for all $0\le k\le \ell-1$. 
Then the next iteration step provides
\begin{align}\label{stabty-ineq}
   \Big| \sum\limits_{i \, j} w_{i,j}\left(h^{i,j} *_{\bb d _{i,j}} v\right)_{\bm n _{\ell}} + y_{\bm n _{\ell}}\Big| &\le    \sum\limits_{i \, j}  w_{i,j} \cdot \noone{h^{i,j}} \cdot \sup\limits_{s < \ell-1 } \lbm v_{\bm n_{s} }\rbm   + \nofty{y}\\ \nonumber
   &\le  \sum\limits_{i \, j} w_{i,j}  \noone{h^{i,j}}  + \nofty{y} \le 2,
\end{align}
where the second inequality  follows from the induction hypothesis. Since for any number $|a|\le 2$, the expression $a-\sign(a)$ lies in the interval $[-1,1]$,  the estimation \eqref{stabty-ineq} shows that $|v_{\bm n_{\ell}}|\le 1$ which finishes the proof by induction. 
\end{proof}

\begin{rem} The stability of the 1st-order weighted $\SD$ schemes \eqref{1st-order-QS-weih_1}-\eqref{1st-order-QS-weih_2} follows from Proposition~\ref{Stability-prop} and the fact that for all 1st-order feedback filters have the form $h^{i,j}\!=\![0, 1]$. 

\end{rem}

We now prove Theorem~\ref{thm-2d-error-bounds-1st-WSD} and Theorem~\ref{thm-2d-error-high-order-WSD}. 
The proof of these results requires the following auxiliary concept and its particular properties. 

 For a given direction $\bm d$,  the generalized directional convolution $h\oast_{\bm d}\Psi$ between a feedback filter $h$ and a function $\Psi$ is defined as 
\begin{equation}\label{gen-convol}
h\oast_{\bm d}\Psi(\bm a)=\sum_{s=1}^L h_s \Psi\big( \bm a -s \bb d\big),  \quad \quad \text{for} \; \bm a \in \R^2. 
\end{equation}
The next lemma suggests  that the generalized convolution $h\oast_{\bm d}\Psi$ can be well approximated by the function derivatives of higher orders once the filter fulfills the moment conditions \eqref{moment-cond}. 

\begin{lem}\label{lem-mon-cond}
Let a feedback filter $h$ be such that $h_n\!=\!0$ for all $n<0$ and $n>L$,  $L\in \N$, and the moment conditions \eqref{moment-cond} are fulfilled. Then for any bivariate  $r$-times differentiable function $\Psi$ and points $\bm a, \bb d\in \R^2$, one has
\begin{equation*}
  h\oast_{\bm d}\Psi(\bm a)  =\Psi( \bm a)+(-1)^{r} C_{h} \sum\limits_{|\alpha|=r}\tfrac{ \bm d^{\alpha}}{\alpha!}\, \pd^{\alpha}\Psi (\bm a)+  \sum_{j=1}^L h_j R_{\bm a, r} \big({-j\bb d}\big), 
\end{equation*}
where  $C_h$ is the filter constant and $R_{\bm a, r}$ is the remainder term in the $r$th order Taylor expansion of  $\Psi$. 
\end{lem}

 \begin{proof} Taylor's formula applied to  $\Psi$ around  $\bm a$ with increment $-j\bb d$ yields
\begin{equation*}
    \Psi(\bm a -j\bb d)=\sum_{\ell=0}^r \sum\limits_{|\al|=\ell}\tfrac{(-j\cdot \bb d )^{\al}}{\al!} \, \da {\al} \Psi(\bm a)+ R_{\bm a, r} \big({-j\bb d}\big), 
\end{equation*}
and, combined with  the definition of the generalized directional convolution~\eqref{gen-convol}, 
\begin{align*}
 h\oast_{\bm d}\Psi(\bm a) 
  & =   \sum\limits_{\ell=0}^r\sum\limits_{|\alpha|=\ell}\tfrac{(-1)^{\ell}   \bb d^{\alpha}}{\alpha!}\, \pd^{\alpha}\Psi (\bm a) \sum_{j=1}^L h_j \cdot j^\ell + \sum_{j=1}^L h_j R_{\bm a, r} \big({-j\bb d}\big).
 \end{align*}
 Since $\delta^0-h$ satisfies the moment conditions \eqref{moment-cond},  the first sum equals $\Psi(\bm a)$ for $\ell=0$ and vanishes for $0<\ell<r$. Rewriting the summand for $\ell=r$ in terms of the filter constant $C_h$,  we obtain the desired result. 
\end{proof}



\begin{proof}[Proof of Theorem~\ref{thm-2d-error-high-order-WSD}] 
We need to bound the supremum norm of 
 the error signal $f_{\lam} - f_q$, which  due to the definition of the quantization scheme  \eqref{rth-order-sim-QR-filter_1}-\eqref{rth-order-sim-QR-filter_2}  can represented as 
\begin{align}
   f_{\lam} (\x) &- f_q (\x)  =  \frac{1}{\lam^2} \sum\limits_{\bm n \in \N^2} \lb y_{\bm n}-q_{\bm n}\rb \Phi\lb \x - \frac{\bm n}{\lam}\rb \\
    & = \frac{1}{\lam^2} \sum\limits_{\bm n \in \N^2} v_{\bm n} \Big( \Phi\lb \x - \frac{\bm n}{\lam}\rb -\sum\limits_{i\,j} w_{i,j} \sum_{s=1}^L h^{i,j}_s \Phi\Big( \x - \frac{\bm n +s \bb d_{i,j}}{\lam}\Big) \Big).    \label{e-t-1-2nd}
\end{align}
For each sum $ \sum_{s=1}^L h^{i,j}_s \Phi\lb \x - \frac{\bm n}{\lam} -\frac{s \bb d_{i,j}}{\lam}\rb$, which is the generalized convolution $h^{i,j} \oast_{\bb d_{i,j}/\lam}\!\Phi$ at  $\bm a _{\bm n}= \x-\frac{\bm n}{\lam}$, Lemma~\ref{lem-mon-cond} provides 
\begin{equation*}
{h^{i,j} \oast_{\frac{ \bb d_{i,j}}{\lam}}\!\Phi(\bm a_{\bm n})}\!=\!\Phi(\bm a_{\bm n})\,\!+\,C_{h^{i,j}}\!\sum\limits_{|\alpha|=r}\tfrac{    \lb -\bb d_{i,j}\rb^{\alpha}}{\alpha!\,\lam^r}\, \pd^{\alpha}\Phi (\bm a_{\bm n})\,\!+\!\,\sum_{s=1}^L h^{i,j}_s R_{\bm a _{\bm n}, r} \big({-\tfrac{s}{\lam}\bb d_{i,j}}\big).
\end{equation*}
 Using this along with the fact that weight matrix elements add up to one, we can represent the quantization error as follows
   \begin{align}
          f_{\lam} (\x) - f_q (\x) &= \frac{1}{\lam^2} \sum\limits_{\bm n \in \N^2} v_{\bm n} \Big( \sum\limits_{i \, j} w_{i,j}(-1)^{r+1}  \sum\limits_{|\alpha|=r}\tfrac{ C_{h^{i,j}} \lb \bb d_{i,j}\rb^{\alpha}}{\alpha! \, \lam^r}\,\pd^{\alpha}\Phi(\bm a_{\bm n}) \nonumber\\
         &\quad +  \sum\limits_{i \, j} w_{i,j} \sum_{s=1}^L h^{i,j}_s R_{\bm a_{\bm n}, r} \big({-\tfrac{s}{\lam}\bb d_{i,j}}\big) \Big) \nonumber\\
          &=\frac{1}{\lam^r} \sum\limits_{|\alpha|=r} \sum\limits_{\bm n \in \N^2} \frac{v_{\bm n}}{\alpha ! \, \lam^2 }\, \pd^{\alpha}\Phi (\bm a_{\bm n}) \cdot \sum\limits_{i \, j} w_{i,j} C_{h^{i,j}}\bb d_{i,j}^{\alpha} \label{err-1-1}\\
          & \quad + \frac{1}{\lam^2}  \sum\limits_{\bm n \in \N^2} {v_{\bm n}} \sum\limits_{i \, j} w_{i,j} \sum_{s=1}^L h^{i,j}_s\, R_{\bm a_{\bm n}, r} \big({-\tfrac{s}{\lam}\bb d_{i,j}}\big).   \label{err-1-2}
   \end{align}
To estimate the product \eqref{err-1-1}, we need to bound the corresponding factors. For the first factor, we obtain that for an appropriate constant $C>0$
\begin{align*}
  \sum\limits_{|\alpha|=r}\bigg( \sum\limits_{\bm n \in \N^2} \frac{v_{\bm n}}{\alpha! \, \lam^2} \, \pd^{\alpha}\Phi \lb \bm a_{\bm n} \rb \bigg)^2 &\le \nofty{v}^2\sum\limits_{|\alpha|=r}\frac{1}{\alpha!}\bigg(  \sum\limits_{\bm n \in \N^2} \frac{|\pd^{\alpha}\Phi \lb \bm a_{\bm n} \rb|}{ \lam^2}  \bigg)^2 \\
  &\le C \nofty{v}^2 \big(\left\|  \nabla^r \Phi\right\|_{1,2}\big)^2 
\end{align*}
where the last inequality follows from the observation that the inner sums under consideration are exactly Riemann sums approximating the integrals that define $\left\|  \nabla^r \Phi\right\|_{1,2}$. The second factor can be written  in terms of $\bb d_{i,j}$ components,
  \begin{equation*}
    \sum_{|\alpha|=r}\Big(  \sum_{i,j}C_{h^{i,j}}w_{i,j}\lb  \bb d_{i,j}\rb^\alpha\Big)^2= \sum_{m=0}^r\Big(  \sum_{i,j}C_{h^{i,j}}\cdot w_{i,j} \cdot i^{r-m}j^m\Big)^2= (C_{\W})^2. 
  \end{equation*}

Exploiting these estimates with Cauchy–Schwartz argument yields  
\begin{equation*}
    \lbm f_{\lam} (\x) - f_q (\x) \rbm \le \frac{1}{\lam^r}  \nofty{v}\Big(  C_{\W}\cdot C \cdot \left\|  \nabla^r \Phi\right\|_{1,2}+ \mathcal{O}(\lam^{-1})\Big), 
\end{equation*}
where the term $\bO(\lam^{-1})$ arises from the integral form of the remainder in Taylor's formula, see Appendix~\ref{Taylor-Rem-Est-1st-order} for details. 
This  completes the proof as this upper bound is independent of $\bm x$.  
\end{proof} 

\begin{rem} Theorem~\ref{thm-2d-error-bounds-1st-WSD} is a direct corollary of Theorem~\ref{thm-2d-error-high-order-WSD} if we set $r\!=\!1$ and $h=[0,1]$, as for this filter one has $C_h=1$.
\end{rem}

  \section{Numerical Experiments}\label{numerics-sec}
  
 The goal of this section it to empirically study the performance of the $2D$ quantization and  error diffusion algorithms motivated by the theoretical investigations in the previous sections.  We first explore the scenario of quantizing bandlimited signals, directly corresponding to our theoretical results, before exploring the performance for digital halftoning the main motivating application of this paper.  
 
  \subsection{Quantization of Bandlimited Functions in $2D$}\label{numerics-sec-fucntions}

For quantization of  bivariate bandlimited functions we summarize our method in  Algorithm~\ref{alg:1}. In our numerical implementations, we use the sinc-function as a kernel even though we are aware that technically it does not satisfy the assumption on  Schwartz  function class. The reason is that  it is simpler to implement and can be closely approximated by Schwartz functions, so we expect comparable behaviour on signals represented by finitely many samples. 

As a test case, we consider the bandlimited  function 
   \begin{equation}\label{test-func}
       f(x_1,x_2)= \frac{3}{10}\cdot \mathrm{Re}\left[\e^{-\im(3x_1+2x_2)}\cos\lb\frac{x_2}{3}\rb\right], 
   \end{equation}
with  frequency support  
in the square $\S_{\frac{10}{3\pi}}$. 
The reconstruction kernel is scaled as follows 
\begin{equation*}
    \Phi(x_1,x_2)=\mathrm{sinc2d}(x_1,x_2):= 25\cdot \mathrm{sinc}\lb 5x_1\rb \cdot\mathrm{sinc}\lb 5x_2\rb.
\end{equation*}
$f$ is sampled on the non-negative numerical lattice 
\begin{equation}
    \L_{\lam}:=\lbc \lb\frac{n_1}{\lam}, \frac{n_2}{\lam} \rb\in \R_+^2\colon n_1,n_2=0, \dots, 10\lam \rbc
\end{equation}
with oversampling $\lam$ rate varying in the set $ \lbc 75+25n, n=0,...,8 \rbc$.
The sample sequence  is denoted by  $y=\lbc y_{n_1,n_2 }\rbc_{n_1,n_2=0} ^{10\lam,10\lam}$ with ${y_{n_1,n_2}:{=}f\lb\frac{n_1}{\lam}, \frac{n_2}{\lam}\rb}$.

\bigskip

\begin{algorithm}[H]

{\small
\KwData{\begin{itemize}
 \item samples  $\quad \quad \quad \quad y=\lbc y_{n_1,n_2}\rbc_{n_1,n_2=0}^{N_1,N_2}$ \quad  with $y_{n_1,n_2}: =f(\frac{n_1}{\lam},\frac{n_2}{\lam})$
    
        \item bounded region $ \quad \mathcal{D}\subset \R_+^2$
\end{itemize}


{\bf Quantization setup:}  
\begin{itemize}
         \item weight matrix  $  \quad \quad \W \in \R^{(\ell+s+1) \times (p+1)} $
        \item  feedback filters $ \quad \quad h^{i,j}\in \ell^{1}(\Z)$
        \item kernel  $  \quad \quad  \quad \quad \quad \quad\Phi \in C^{\infty}(\R^2)$ \quad  with $
    \F\Phi(\bm \xi)=\left\{\begin{matrix}
    1,& \bm \xi \in S_\lam,\\
    0,&\bm \xi \notin S_\lam.\end{matrix}\right.
$
\end{itemize}

 \vspace{-5mm}
}

\vspace{2mm}

\Begin{

\vspace{1mm}

 for $ n_1= 0\dots N_1$, \\
  \quad\quad  $n_1=1,..., N_2$ \\
    
     \vspace{-5mm}
    
        \begin{align*}
  v_{n_1,n_2}&-\sum\limits_{i \, j}w_{i,j}\left(h^{i,j}*_{\bb d_{i,j}} v\right)_{n_1,n_2} =y_{n_1,n_2}- q_{n_1,n_2}\\[-2pt]
 q_{n_1,n_2}&= \sign\Big(   \sum\limits_{i \, j} w_{i,j}\left(h^{i,j}*_{\bb d_{i,j}}v\right)_{n_1,n_2} + y_{n_1,n_2}\Big) 
 \end{align*}
  \vspace{-5mm}
 }
 
\KwResult{\begin{itemize}
        \item 1-bit samples  $   \; \quad \quad \quad \quad \quad q=\lbc q_{n_1,n_2}\rbc_{n_1,n_2=0}^{N_1,N_2}$
        \item  quantized rep.   $ \quad \quad \quad \quad \quad \quad \;\ds  f_q(\x) $
\item  quantization error   $ \quad \quad \quad  \ds   err:= \max\limits_{\x\in \mathcal{D}_{dics} }|f_{\lam}(\x) - f_q(\x)| $
\end{itemize}
}

\label{alg:1}
 \caption{\it \small Weighted~$\SD$~Quantization Schemes for $2D$ Functions}
}
\end{algorithm}

  \bigskip

\begin{table}[h!]
    \centering
    \begin{tabular}{|c|p{3cm}|c| c|}
    \hline
    Quantization scheme  & Weight matrix $\W$ &  Constant $C_{\bf W}$  & Abbrv.
    \\[0.5ex] 
 \hline\hline
     Row-by-row    &   $\begingroup 
\setlength\arraycolsep{2pt}
\begin{pmatrix}
     \bm{ 0} & 1\\
     0& 0  \end{pmatrix}
\endgroup $ & $1$& $1st$-$RbR$ \\[2ex] 
 \hline
   Floyd–Steinberg   &   $ \begingroup 
\setlength\arraycolsep{2pt}
\begin{pmatrix}
     0 & \bm 0 & \frac{7}{16}\\[1.5pt] 
     \frac{3}{16} & \frac{5}{16}& \frac{1}{16}  \end{pmatrix}
\endgroup $ & $\frac{\sqrt{106}}{16}\approx0.65$ &$F$-$S$\\[2ex] 
 \hline
Shiau-Fan    &
 $ \begingroup 
\setlength\arraycolsep{2pt} \begin{pmatrix}
  0& 0 & 0 & \bm{ 0} & \frac{8}{16}\\[1pt] 
\frac{1}{16} &\frac{1}{16}     & \frac{2}{16} & \frac{4}{16}& 0 \end{pmatrix} \endgroup 
$ & $\frac{\sqrt{65}}{16}\approx0.5$ 
&$S$-$Fan$\\[2ex] 
 \hline
     Averaged & $ \begingroup 
\setlength\arraycolsep{2pt}
\begin{pmatrix}
      \bm{ 0} & \frac{1}{2}\\
      \frac{1}{2}& 0 \end{pmatrix}
\endgroup $ & $\frac{1}{\sqrt{2}}\approx0.71$&$1st$-$A$\\[2ex] 
 \hline
      Optimal-2&  $ \begingroup 
\setlength\arraycolsep{2pt}
\begin{pmatrix}
   0& \bm{ 0} & \frac{7}{10}\\
   \frac{3}{10} & 0&0  \end{pmatrix}
\endgroup $ &$\frac{1}{\sqrt{5}}\approx0.45$ &$Opt$-$2$\\[2ex] 
 \hline
   Optimal-4& $ \begingroup 
\setlength\arraycolsep{2pt}
\begin{pmatrix}
  0& 0&0 &0&  \bm{ 0} & \frac{21}{26}\\
   \frac{5}{26}&0& 0&0& 0& 0  \end{pmatrix}
\endgroup $ &$\frac{1}{\sqrt{26}}\approx0.19$&$Opt$-$4$\\
   \hline
    \end{tabular}
    
        \vspace{-2mm}
        
    \caption{ \it First-order weighted $\SD$ quantization schemes. The element $w_{0,0}$ is denoted by zero in bold. }
    \label{tab1}
\end{table}


\begin{table}[h!]
    \centering
    \begin{tabular}{|c|c|c|c| c|}
    \hline
    \makecell{  QS,\\
    Abbrv.} &  $\W$ & Filters $h^{i,j}$& Ext.\,weigh.\,matrix& \makecell{$C_{\bf W}$\\ $\approx$}
    \\[0.5ex] 
 \hline\hline
  { \small \makecell{  Row-by-row \\
 $2nd$-$RbR$}  }  &   $ \begingroup 
\setlength\arraycolsep{2pt}
\begin{pmatrix}
     \bm{ 0} & 1\\
     0& 0  \end{pmatrix}
\endgroup $ &  {\small ${ h^{0,1}{=}[0, 4/3, 0, 0, -1/3]}$}  & $ \begingroup 
\setlength\arraycolsep{2pt}
\begin{pmatrix}
     \bm{ 1} & -\frac{4}{3} & 0 & 0& \frac{1}{3}\\
            0& 0 & 0 & 0 & 0 \end{pmatrix}
\endgroup $  &  $1,68$
\\[2ex] 
 \hline
    { \small   \makecell{ Averaged-33\\ \\ $2nd$-$A33$ }} & $ \begingroup 
\setlength\arraycolsep{2pt}
\begin{pmatrix}
     \bm 0 & \frac{1}{2}\\
      \frac{1}{2}& 0 \end{pmatrix}
\endgroup $ & {\small \makecell{  ${h^{1,0}{=}h^{0,1}}$ \\ ${ h^{0,1}{=}[0, 4/3, 0, 0, -1/3]}$} } 
& $\begingroup 
\setlength\arraycolsep{2pt}
\begin{pmatrix} 
     \bm{ 1} & -\frac{4}{6} & 0 & 0& \frac{1}{6}\\[-1pt]
     -\frac{4}{6}&0& 0&0 &0\\ 0 & 0& 0&0 &0 \\[-2pt]
     0 & 0& 0&0 &0\\[-2pt]
\frac{1}{6} &0& 0&0 &0\\\end{pmatrix}
\endgroup
$& $1,17$ 
\\[2ex] 
\hline
{\small  \makecell{Averaged-34\\ \\ $2nd$-$A34$ }} & $ \begingroup 
\setlength\arraycolsep{2pt}
\begin{pmatrix}
     \bm 0 & \frac{1}{2}\\
      \frac{1}{2}& 0 \end{pmatrix}
\endgroup $ &{\small \makecell{   ${ h^{0,1}{=}[0, 4/3, 0, 0, -1/3]}$\\
$h^{1,0}{=}\quad \quad  \quad \quad \quad\quad  \quad  \quad$\\$ \quad [0, 5/4, 0, 0, 0, -1/4]$}  }
& $\begingroup 
\setlength\arraycolsep{2pt}
\begin{pmatrix} 
     \bm{ 1} & -\frac{4}{6} & 0 & 0& \frac{1}{6}\\[-1pt]
     -\frac{5}{8}&0& 0&0 &0\\ 0 & 0& 0&0 &0 \\[-1pt]
     0 & 0& 0&0 &0\\[-2pt]
      0 & 0& 0&0 &0\\[-2pt]
\frac{1}{8} &0& 0&0 &0\\\end{pmatrix}
\endgroup
$ &$1,27$ \\[2ex] 
   \hline
    \end{tabular}
    
        \vspace{-2mm}
    
    \caption{ \it Second-order weighted $\SD$ quantization schemes. The element $w_{0,0}$ is denoted in bold. }
    \label{tab2}
\end{table}

\newpage 

To measure the resulting quantization error, we compare,  for values $(x_1, x_2)$ in the square ${[2, 8 ]\times[2, 8 ]}$, the values of the approximant    
\begin{equation}\label{f_approx_num_exp}
    f_{\lam}(x_1,x_2):=  \frac{1}{ \lam^2}\sum\limits_{ n_1=0}^{10\lam}\sum\limits_{ n_2=0}^{10\lam} y_{n_1,n_2}\cdot\,  \mathrm{sinc2d}\Big(x_1-\frac{n_1}{\lam},x_2-\frac{n_2}{\lam} \Big),
\end{equation}
and the quantized representative $f_{q}$, 
which is computed by substituting in \eqref{f_approx_num_exp} $y_{n_1,n_2}$ by the corresponding 1-bit sample $q_{n_1,n_2}$, for all $n_1,n_2\in \{0,\dots,10\lam\}$.

 In our numerical experiments, the 1-bit samples $q\!=\!\lbc q_{n_1,n_2 }\rbc_{n_1,n_2=1} ^{10\lam,10\lam}$ of $f$ are constructed  using weighted $\SD$ schemes compiled  in Table~\ref{tab1} (first order) and Table~\ref{tab2} (second order). 
 The Optimal-2 and Optimal-4 schemes in Table~\ref{tab1} are the weighted $\SD$ scheme resulting from the solution of Optimization Problem \ref{opt-prob-1st-or} 
 with size parameter $s\!=\!2$ or $s\!=\!4$, respectively.

\begin{figure}[ht!]
 \centering
 
  \hspace*{\fill}%
  
  \subcaptionbox{$f$ and its approximant $f_{\lam}$;\\ approx.\,error ${4,848\cdot10^{-3}}$.     \label{fig32:b}}{\includegraphics[width=2.3in]{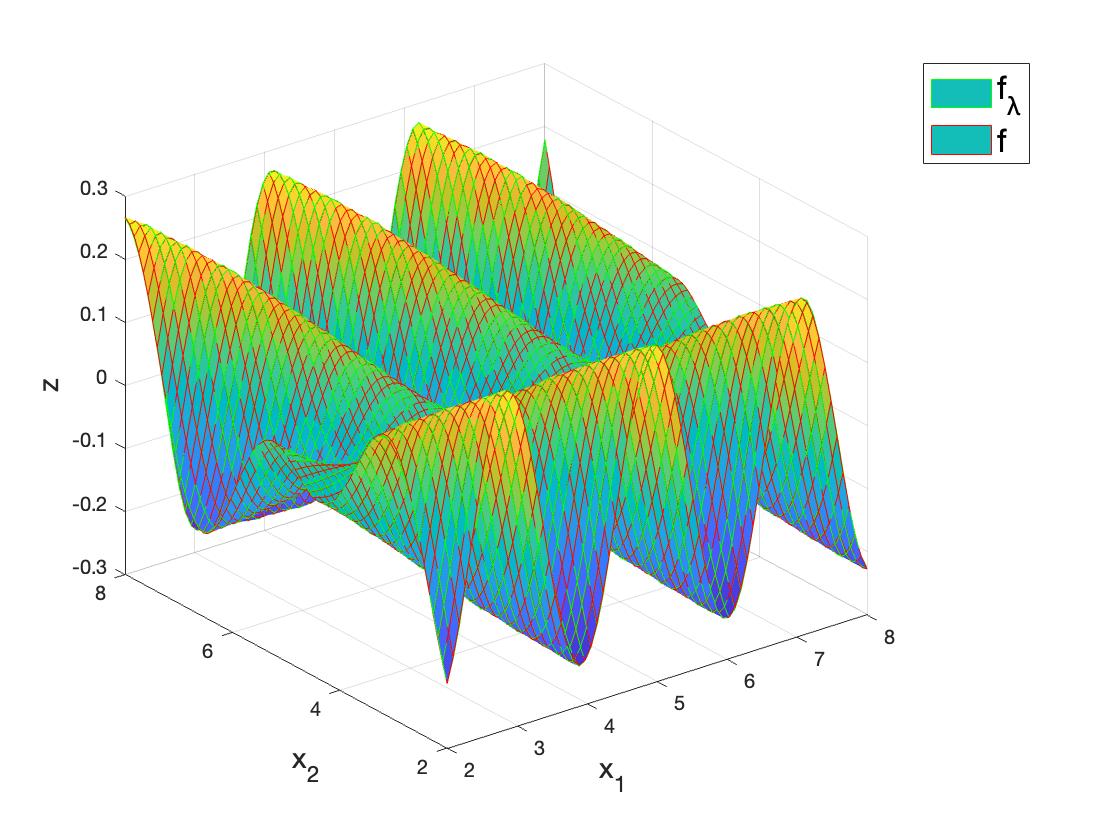}}%
  \hfill%
  \subcaptionbox{ Error signal 
  for  $1st$-$RbR$ scheme;  \\ maximal  amplitude  ${2,251\cdot10^{-2}}$. \label{fig33:a}}{\includegraphics[width=2.3in]{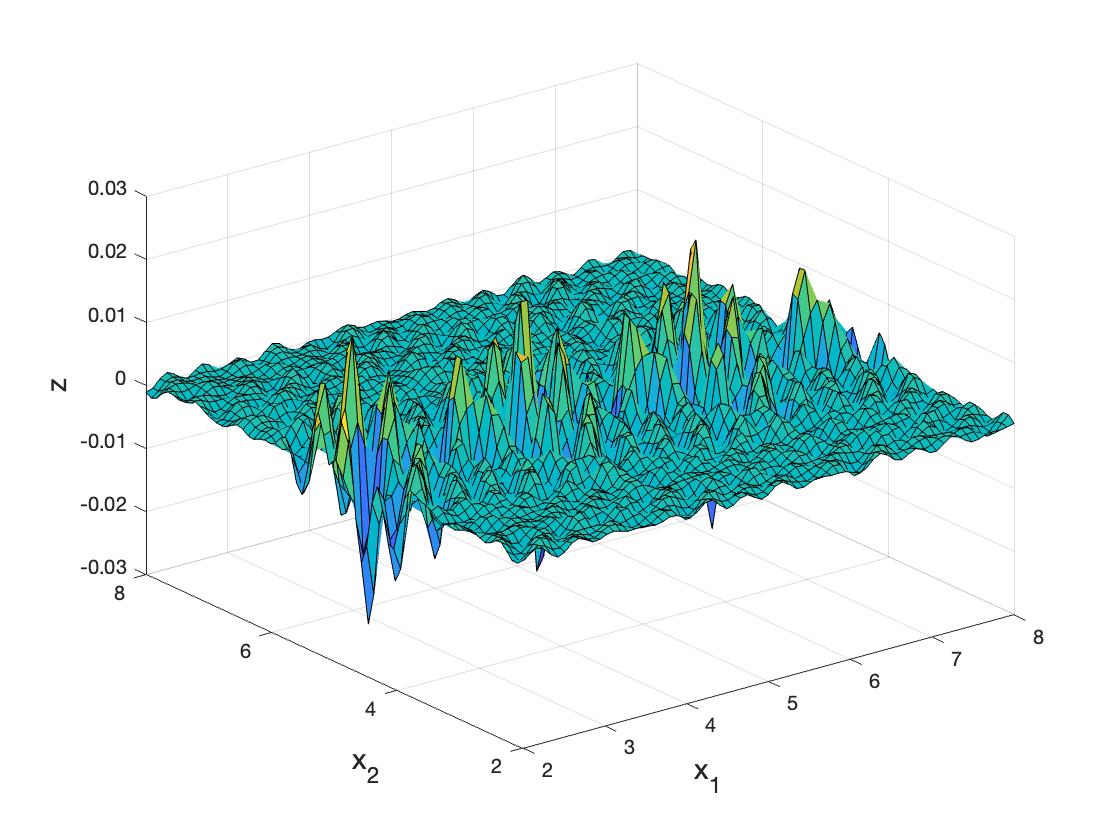}}
  
 \centering
 
   
    \subcaptionbox{ Error signal 
    for  {$1st$-$A$} scheme; maximal amplitude  ${1.293\cdot10^{-2}}$.    \label{fig34:b}}{\includegraphics[width=2.3in]{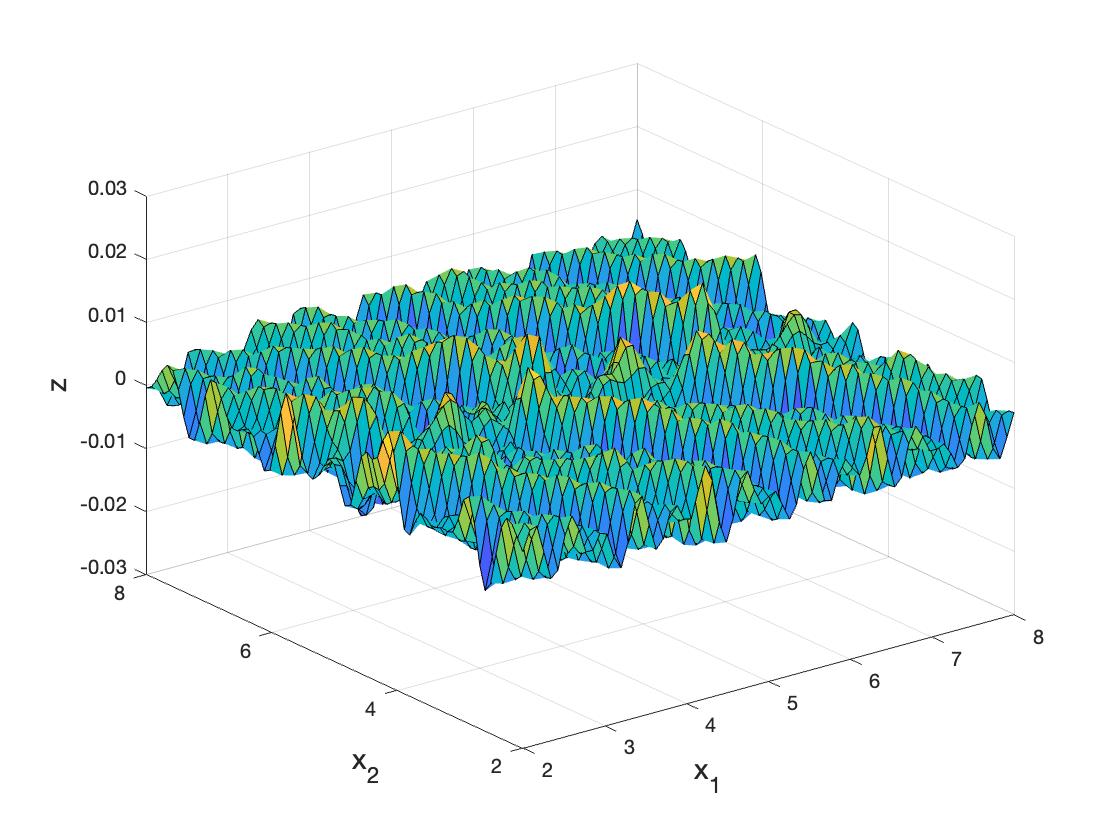}}
  \hfill%
   \subcaptionbox{ Error signal for $Opt$-$4$ scheme; maximal  amplitude  ${4,663\cdot10^{-3}}$.    \label{fig34:a}}{\includegraphics[width=2.3in]{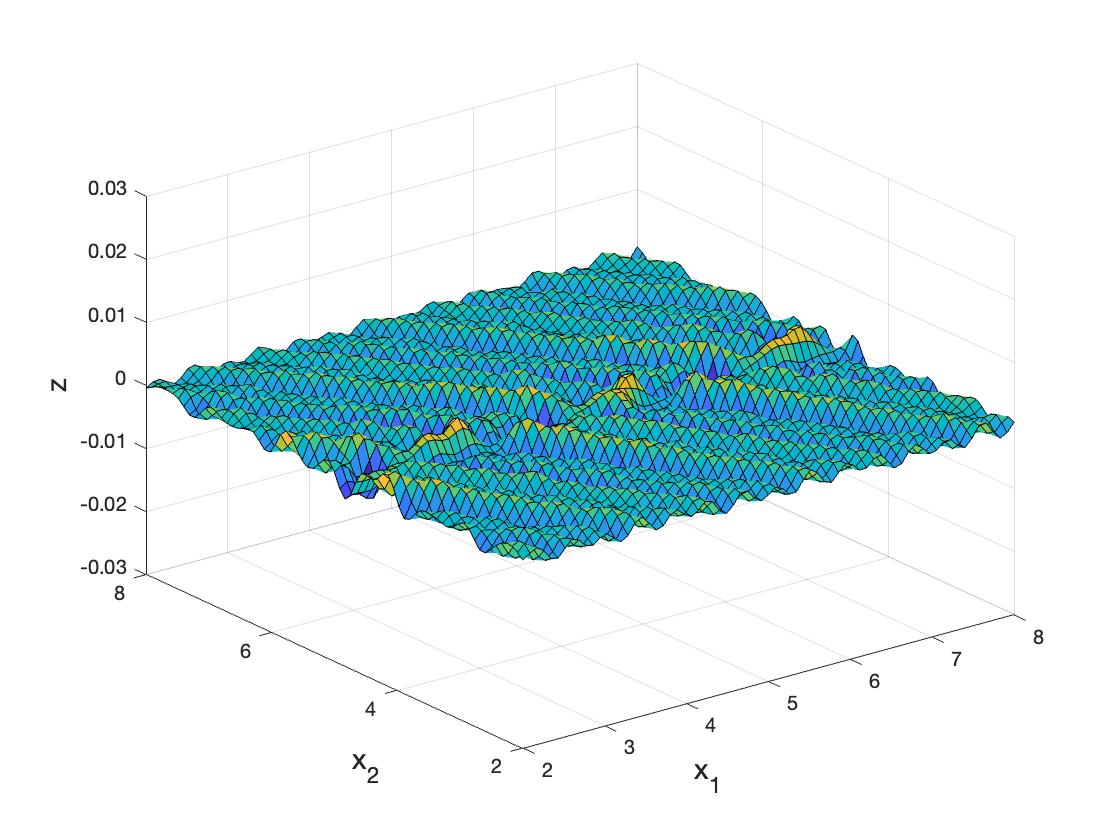}}
 %


    \vspace{-2mm}

  \caption{ \it Performance of 1st-order weighted $\SD$ quantization schemes 
  with oversampling rate
  ${\lam=150}$ for a bandlimited signals. The experiment demonstrates the  benefits of optimizing the weights in this context.}
  
 \label{fig:num:3}
\end{figure}

Figure \ref{fig:num:3} illustrates the performance of  Algorithm~\ref{alg:1} with a first order feedback filter for three different weight matrices and  oversampling rate of $\lam\!=\!150$ applied to  $f$ as given in  \eqref{test-func}. 
Setting $\lam\!=\!150$ produces the approximant $f_{{\lam}}$ representing $f$ with maximal error of order $10^{-3}$, see Figure~\ref{fig32:b}. 
We observe that applying one-dimensional schemes row-by-row gives rise to difficulties in the area with actively varying function values, which illustrates the advantage of weighting different directions. At the same time, we observe an additional advantage for the weight minimizing the constant $C_{\W}$.

\begin{figure}[ht!]

  \subcaptionbox{ First-order schemes \label{fig51:a}}{\includegraphics[width=2.3in]{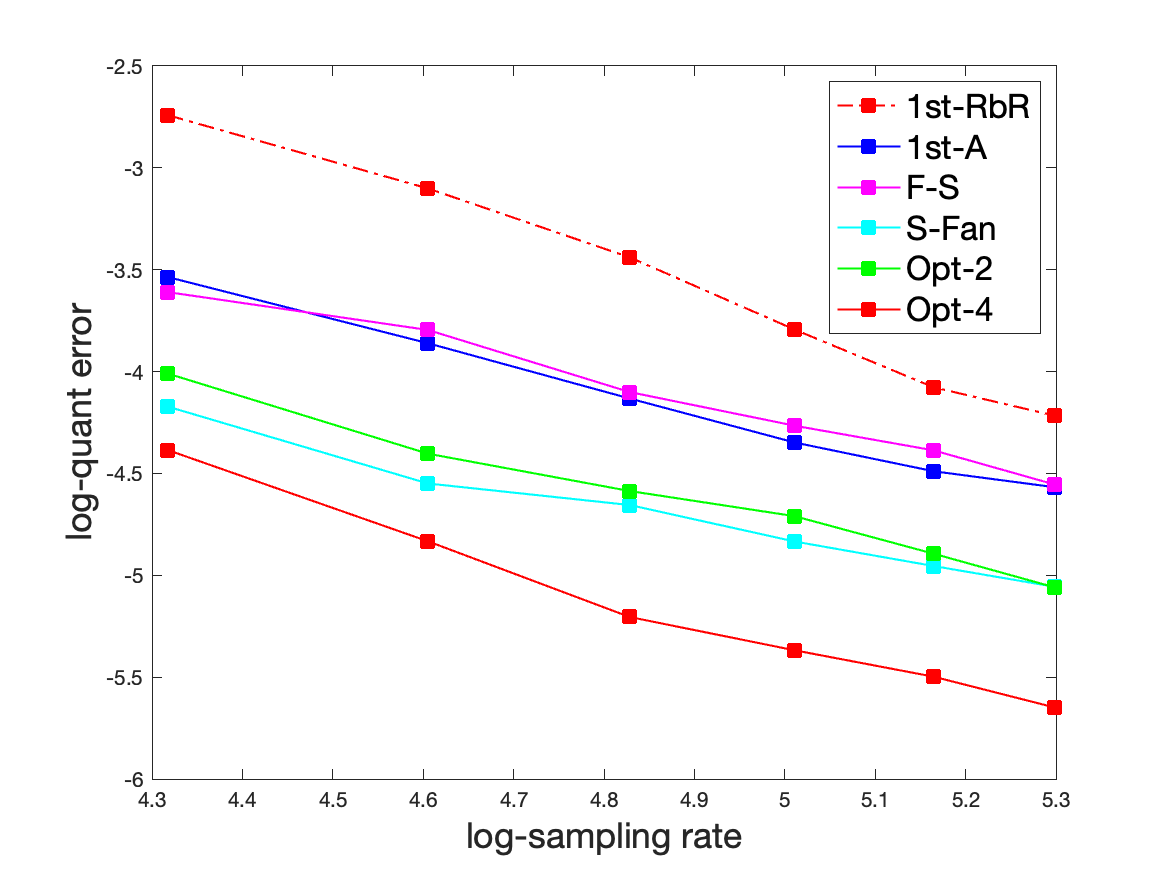}}
  \hfill%
  \subcaptionbox{  Second-order schemes.         \label{fig52:b}}{\includegraphics[width=2.3in]{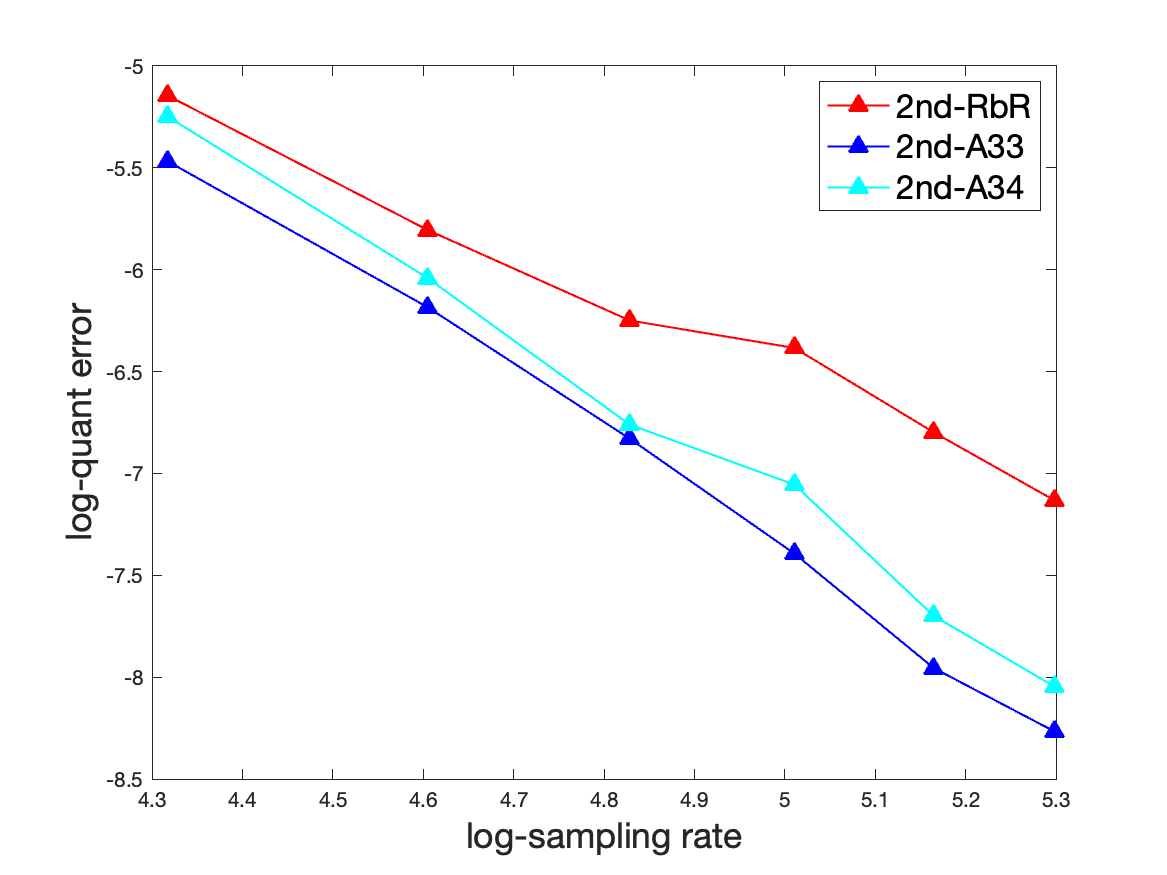}}

  \vspace{-2mm}
  
  \caption{ \it Maximal  quantization error in domain $[2,\, 8]\times [2,\, 8]$ for weighted $\SD$ quantization schemes with oversampling rate $ {\lam \in \lbc 75+25n, n=0,...,8 \rbc}$.}
  
 \label{fig:num:4}
\end{figure}

We explore this observation in more detail in a line of experiments summarized in 
 Figure~\ref{fig:num:4} comparing the performance of 1st- and 2nd-order $\SD$ schemes  with different weight matrices  for various values of the oversampling rate $\lam$. While our findings confirm that the constant $C_{\W}$ identified in our analysis appropriately captures  the performance  differences between 1st-order schemes with different weight matrices, we also find that all first-order schemes are outperformed by any second-order scheme under consideration.

  \subsection{Digital Halftoning of Images and Similarity Indices}\label{num-halftone-sec}
  
As explained above, our observation that many error diffusion schemes can be interpreted as weighted 1st-order $\SD$ schemes,  motivated us to mathematically analyze such schemes and systematically explore their use for digital halftoning. In particular, driven by their superior performance for bandlimited signals, we aim to adapt weighted second-order schemes for this purpose. The $\SD$ perspective then motivates to aim for achieving stability in the sense discussed in the previous section, as otherwise, one expects the error to accumulate. 

In this subsection, we confirm these heuristics  using numerical experiments  for a wide variety of images, demonstrating  that indeed ``stabilized'' weighted  second-order  $\SD$ schemes  
outperform a number of state-of-the-art error diffusion schemes in terms of the visual similarity between the halftoned image and the original.

In our experiments, we consider both color and gray-scale images. We now explain our setup for color images, gray-scale images are treated analogously.

  We represent color images  as RGB matrices, $\mathcal{I}_{RGB}\!:=\!\left\{\mathcal{I}^{RGB}_{n_1,n_2} \right\}_ {n_1,n_2=0}^{N_1,N_2}$, consisting of three color channels, each given as a sample array. In order to construct a halftoned counterpart of $\mathcal{I}_{RGB}$, we propose to use Algorithm~\ref{alg:2}, as introduced below. 
  
 In particular, we start by  converting  the image ${\mathcal{I}_{RGB}}$  to  \texttt{double} ${\mathcal{I}_{dbl}\!\in\![0,1]^{N_1\times N_2\times 3}}$, rescaling it to the range $[-1.15,0.85]$ using the~sharpening~map
  $$
\mathcal{I}:=\texttt{sharp}\big(2\,\mathcal{I}_{dbl}-1.15\big), \quad \text{with} \quad \texttt{sharp}(x)=\left\{\begin{matrix}
    -1,& x \le -1,\\
    1,&x \ge 1.\end{matrix}\right.
$$  
This sharpening step is introduced to improve the color fidelity of the halftoned image. 
Then, for each color channel one constructs a $1$-bit array using a weighted $\SD$ scheme.  

In our experiments, we compare weighted $\SD$ schemes of first order (including some of the error diffusion schemes proposed in the literature), of  second order, and of mixed order, that is, schemes applying $\SD$ quantizers of different orders in different directions. 
   
   While weighted $\SD$ schemes of first order are intrinsically stable, schemes of second order or mixed order are typically not stable unless carefully chosen feedback filters are employed. 

 \bigskip
 
\begin{algorithm}[H]
{\small
\KwData{\begin{itemize}
         \item RGB image  $ \quad \quad \quad \mathcal{I}_{RGB}:=\left\{\mathcal{I}^{RGB}_{n_1,n_2} \right\}_ {n_1,n_2=0}^{N_1,N_2}$ 
     \item conversion    $\;\, \quad \quad \quad \mathcal{I}_{dbl}:= \texttt{im2double}(\mathcal{I}_{RBG}) \in [0,1]^{N_1\times N_2\times 3}$
      \item sharping  \quad \quad \quad \quad  $\mathcal{I}=C\cdot\texttt{sharp}(2\,\mathcal{I}_{dbl}-1.15) \in  [-1,1]^{N_1\times N_2\times 3}$
      
\end{itemize}


{\bf Quantization setup:}  
\begin{itemize}
         \item weight matrix  $  \quad \quad \W \in \R^{(\ell+s+1) \times (p+1)} $
        \item  feedback filters $ \quad \quad \quad h^{i,j}\in \ell^{1}(\Z)$

        \item  stability constant $ \quad \quad{\small C=\left\{\begin{matrix}
    0.999,& \text{for } \; 2nd\text{-}SD,\\
    1,&  \text{for }\text{other schemes}.\end{matrix}\right.}$
\end{itemize}

\vspace{-3mm}
}


\Begin{
 for $c= R,G,B$\\
   \quad\quad  $n_1=1,..., N_1$ \\
  \quad \quad \quad  $n_2=1,..., N_2$
  
  \vspace{-5mm}
  
        \begin{align*}
 \quad \quad \quad  v^c_{n_1,n_2}-&\sum\limits_{i \, j}w_{i,j}\left(h^{i,j}*_{\bb d_{i,j}} v^c\right)_{n_1,n_2} =\mathcal{I}^c_{n_1,n_2}- q^c_{n_1,n_2}\\[-2pt]
 q^c_{n_1,n_2}&= \sign\Big(   \sum\limits_{i \, j} w_{i,j}\left(h^{i,j}*_{\bb d_{i,j}}v^c\right)_{n_1,n_2} + \mathcal{I}^c_{n_1,n_2}\Big) 
 \end{align*}
 
 \vspace{-5mm}
 
 }
 
\KwResult{\begin{itemize}
         \item 1-bit image $ \quad \quad \quad \quad q=\left\{q^c_{n_1,n_2} \right\}_{n_1,n_2=1}^{N_1,N_2}\in \{-1,1\}^{N_1\times N_2\times 3}$ 
        \item $bmp$-image     \quad \quad \quad \quad   $\mathcal{I}_q \in \{0,255\}^{N_1\times N_2 \times 3}$
\item  halftoning error  $ \quad \quad \quad  \ds   err_{im}= FSIM(\mathcal{I},\mathcal{I}_q) $
\end{itemize}
}

 \caption{Weighted $\SD$ Schemes  for Digital Halftoning of  Images}
 \label{alg:2}
}
\end{algorithm}

   \bigskip


In this paper we will work with the family of second-order filters with minimal support  given by 
 \begin{equation}\label{k-tab-filter}
     h^\kappa=(0,h^\kappa_1, 0, \dots, 0, h^\kappa_\kappa), \quad \text{for}  \quad h^\kappa_1=\frac{\kappa+1}{\kappa },\quad   h^\kappa_\kappa=-\frac{1}{\kappa}
 \end{equation}
 with $\noone{h^\kappa}=1+\frac{2}{\kappa }$, as introduced in \cite{Gunturk2003}. For the  choice  of the parameter $\kappa$, there is a trade-off. On the one hand, choosing $\kappa$  large will increase the range of applicability of the stability guarantee in Proposition~\ref{Stability-prop}, on the other hand, very large filters increase boundary effects. 
 
 To resolve this issue, we propose to combine $\SD$  schemes built out of 
 both $h^\kappa$ with large and small $\kappa$, properly weighted so that stability for an input amplitude close to one can reached. Two examples of such weighted $\SD$ schemes are described in Table~\ref{tab3}. For instance, the 2nd-order scheme $2nd\text{-}SD$  combines the  filters $h^{550}$ and $h^3$.  Here, the choice of the filter $h^{550}$ is the result of numerically comparing  $h^\kappa$ with $\kappa \in \{100+50k: \; k\!=\!0,\dots16\}$.  While for the resulting  combined filter, Proposition~\ref{Stability-prop}
guarantees stability for signals of maximal amplitude $0.96$,  we  numerically observe stability for amplitudes up to $0.999$. We accomodate for this limitation by rescaling the images with the factor  $0.999$; this has basically no effect on the perceived image, yet ensures stability.  Also for the initial value for the state variable we explored various options and observed best performance for a uniform distribution on $[-0.9.0.9]$. The randomness helps avoid the occurrence of many zero state variable values in the first quantization steps.
 
 The weighted $\SD$ scheme $S\text{-}Fan\text{-}12$ is a combination of 1st and 2nd-order schemes and designed as a reinforcement of the $1$st-order Shiau-Fan scheme by adding the filter $h^3$ in several directions.  Again, stability  can be guaranteed for input signals of amplitude at most  $0.96$, but in our numerical experiments  we do not encounter any instabilities and hence do not propose a rescaling.

\begin{table}[h!]
    \centering
    \begin{tabular}{|@{}0c<{}|@{}0c@{}|@{}0c@{}|@{}0c@{}|}
    \hline
    \makecell{ \small  Quant.\\
   \small Scheme} &  \small  Weight matrix $\W$ &  $(i,j)$-Index sets  & Filters $h^{i,j}$ 
    \\[0.5ex] 
 \hline
 \hline
  { \small  $2nd\text{-}SD$ } & $\begingroup 
\begin{pmatrix} 
     0 & \bm 0 & \frac{88}{199} &\frac{5.5}{199}\\[2pt] 
   \frac{12}{199}  &\frac{87}{199} &\frac{1}{199}&0\\[2pt] 
     0 & \frac{5.5}{199}& 0&0 \end{pmatrix}
\endgroup$ & \begin{tabular}{@{}0c@{}}  \small$ \quad  (0,1), (1,j)$, $j\!\in\!\{-1,0,1\} \quad $  \\ \hline
 \small$(2,0), (0,2)$\end{tabular} &   \begin{tabular}{>{}0c<{}}  \small  \small $h^{i,j}=h^{550} \,$ \\ \hline \small $h^{i,j}=h^{3}$ \end{tabular}
\\[2ex] 
\hline

{ \small$S\text{-}Fan\text{-}12$} & $\begingroup 
\setlength\arraycolsep{2pt}
\begin{pmatrix} 
 0 &   0 &  0 & \bm 0 & \frac{21}{50} &\frac{3}{100}\\[2pt] 
 \frac{2}{50}  &\frac{2}{50} &  \frac{5}{50}  &\frac{17}{50} &0&0\\[2pt] 
  0 &\frac{0.5}{100}&   \frac{0.5}{100} & \frac{2}{100}& 0&0 \end{pmatrix}
\endgroup$
&  \begin{tabular}{>{}0c<{}} \small  $(0,1), (1,\text{-}j)$, $j\!\in\!\{0,1,2,3\} $ \\ \hline \small $(0,2),(2,\text{-}j)$, $j\!\in\!\{0,1,2\} $ \end{tabular}
%
&  \begin{tabular}{>{}0c<{}}  \small $h^{i,j}=[0, 1]$ \\[-1pt] \hline \small
$h^{i,j}=h^{3}$ \end{tabular} \\
\hline

   \hline
    \end{tabular}
    
    \vspace{-2mm}
    
    \caption{ \it Quantization schemes for digital halftoning. The element $w_{0,0}$ is denoted in bold. }
    \label{tab3}
\end{table}

 For  our performance analysis of weighted $\SD$ techniques for digital halftoning, we use 50 distinct color images of the size $1920\times1280$ and their gray-scale counterparts. To measure the quality of the resulting halftoned images, we compute for each image and its halftoned version the Feature Similarity Index (FSIM)  \cite{Zhang2011}. The values of FSIM range between $0$ and $1$, where  $1$ indicates two identical images, and the more dissimilar two images are, the smaller is the corresponding FSIM. 
 
 \begin{figure}[ht!]
 \centering
 
  \hspace*{\fill}%
  
  \subcaptionbox{ FSIM for gray-scale images.  \label{fig41:a}}{\includegraphics[width=2.35in]{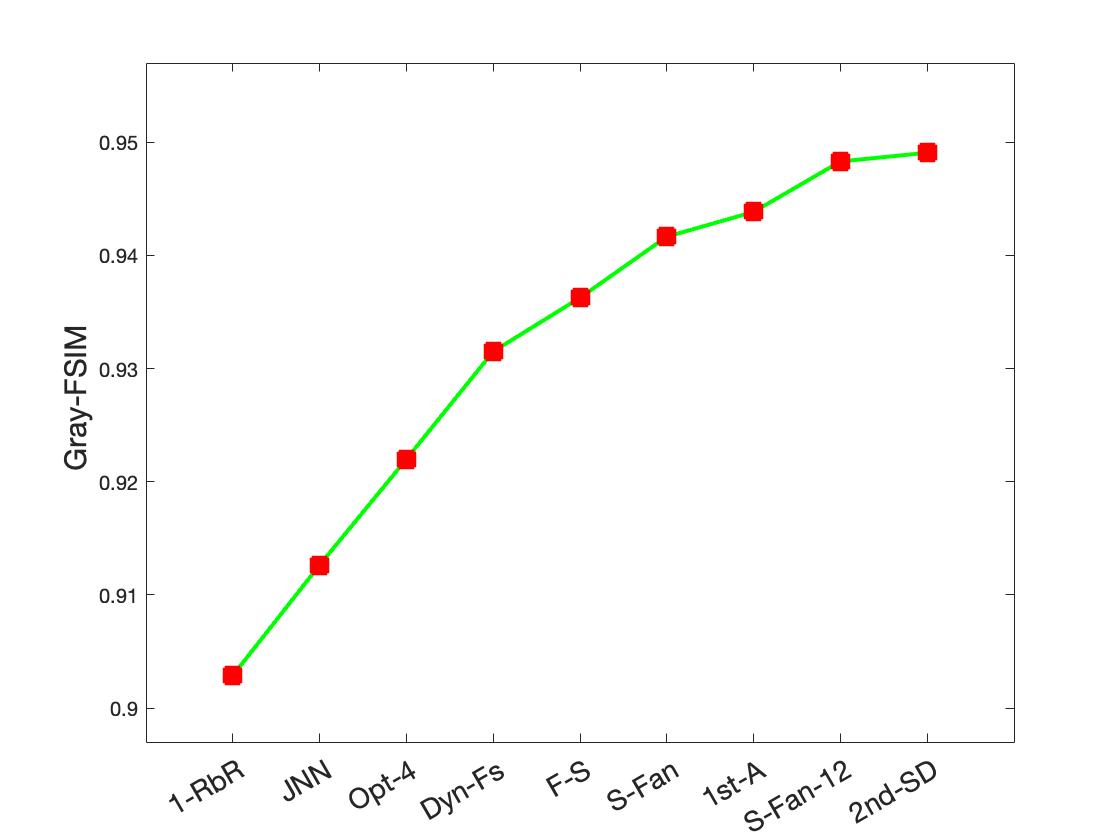}}
  \hfill%
  \subcaptionbox{ FSIM for RGB images.     \label{fig42:b}}{\includegraphics[width=2.35in]{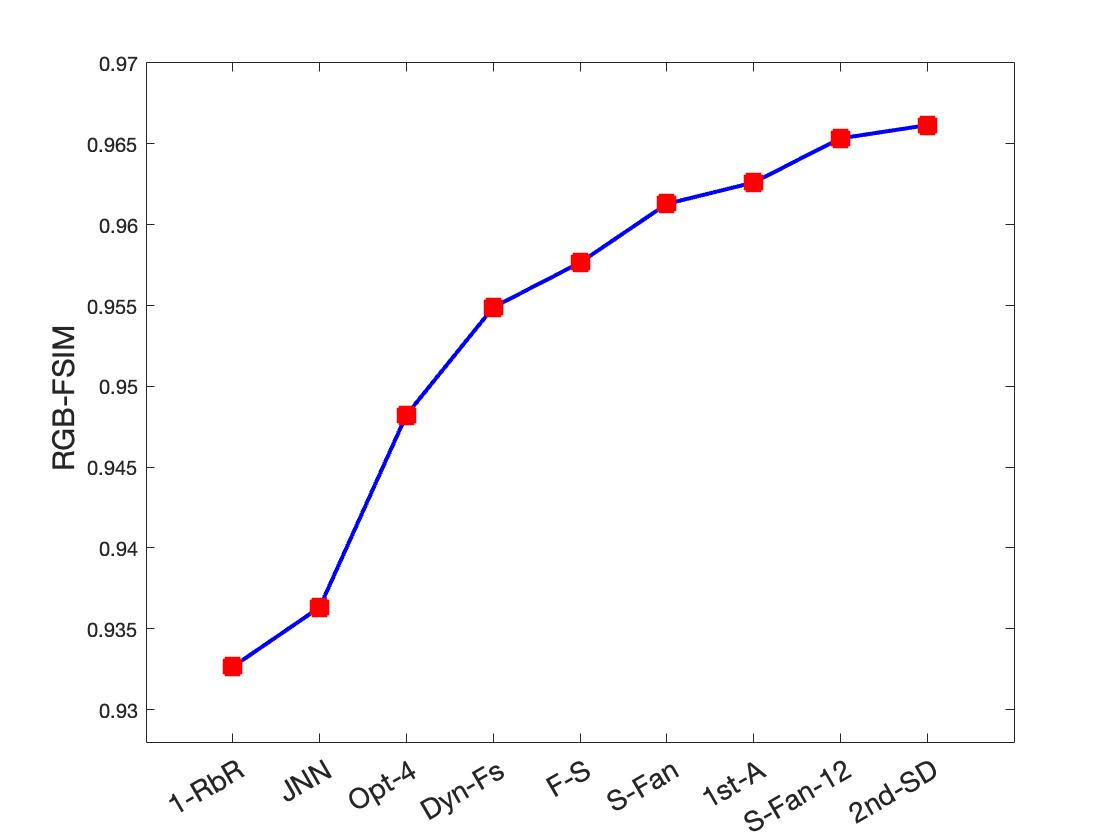}}%
  
  \vspace{-2mm}
 
  \caption{ \small \it  Averaged FSIM over 50 images and their halftoned counterparts generated by state-of-the-art error diffusion schemes as well as the novel schemes proposed in this paper. Here, $JJN$ is the Jarvis-Judice-Ninke scheme  \cite{Jarvis1976}, $Dyn$-$Fs$ is the scheme with dynamic filters,  \cite{Ostromoukhov2001}, and the weighted $\SD$ schemes can be found in Table~\ref{tab1}~and Table~\ref{tab3}. We observe that the best performance in this quality measure is achieved by weighted $\SD$ schemes with second order building blocks as proposed  above.  }
  
 \label{fig:num:8}
\end{figure}    


\begin{figure}[ht!]
\centering
\sbox{\bigpicturebox}{%
  \scalebox{1}[1]{\includegraphics[width=.7\textwidth]{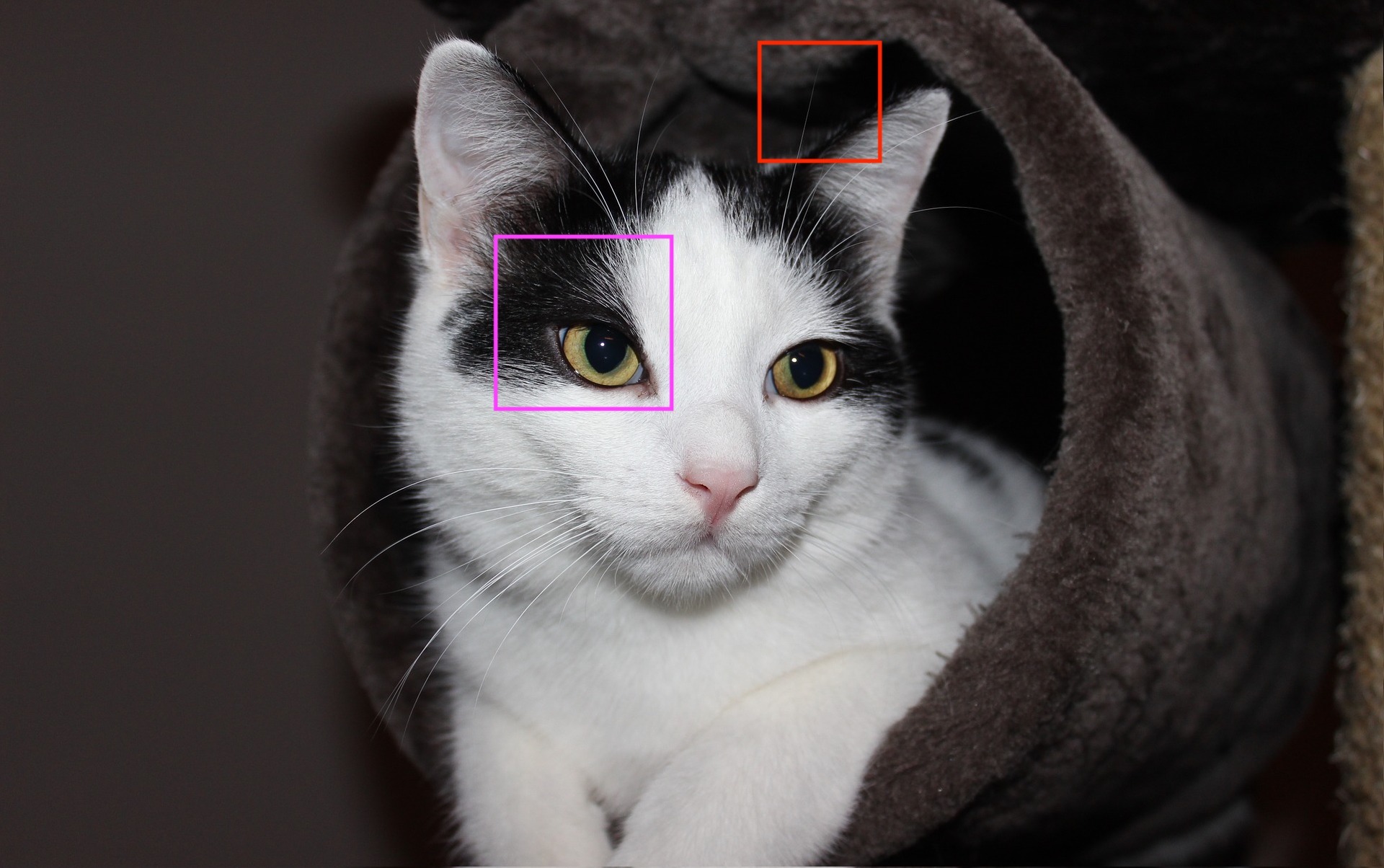}}}
\usebox{\bigpicturebox} \hfill
\begin{minipage}[b][\ht\bigpicturebox][s]{.20\textwidth}
\includegraphics[width=1.1\textwidth]{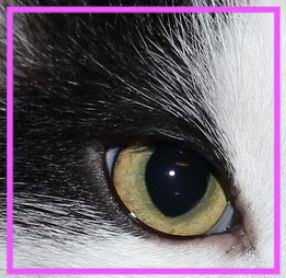}
\vfill
\includegraphics[width=1.1\textwidth]{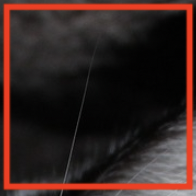}
\end{minipage}
\caption{Original image of size $1920\times1280$ and its corresponding parts zoomed.}
\label{fig:9}
\end{figure}

 Figure~\ref{fig:num:8} shows the average FSIM over 50 images for different $\SD$-based halftoning techniques applied to color and gray-scale images. As one can see, there is almost no difference in the results  obtained for gray-scale and color images. 
 Among the 1st-order techniques, the simple average $1st\text{-}A$ performs the best in terms of FSIM, closely followed by the Shiau-Fan halftoning technique. Images with the smallest similarities are produced by the row-by-row scheme and although,  the  1st-order  $Opt$-$4$ scheme shows exceptionally good performance in terms of the supremum norm for quantizing bandlimited signals, for digital image halftoning it is outperformed by simple averaging. 
 
 Despite minimally rescaling the image amplitude, we encounter the best performance among all weighted $\SD$ scheme for the 2nd-order scheme $2nd\text{-}SD$.  The second best result is produced by the reinforced Shiau-Fan scheme $S\text{-}Fan\text{-}12$ of mixed order. In this case, we see similar behavior of weighted $\SD$ schemes both for quantization of bivariate bandlimited functions and digital halftoning of images, namely, higher-order schemes perform better then the ones of first~order. 

{
\begin{figure}[ht!]
\centering
\sbox{\bigpicturebox}{%
  \scalebox{1}[1]{\includegraphics[width=.7\textwidth]{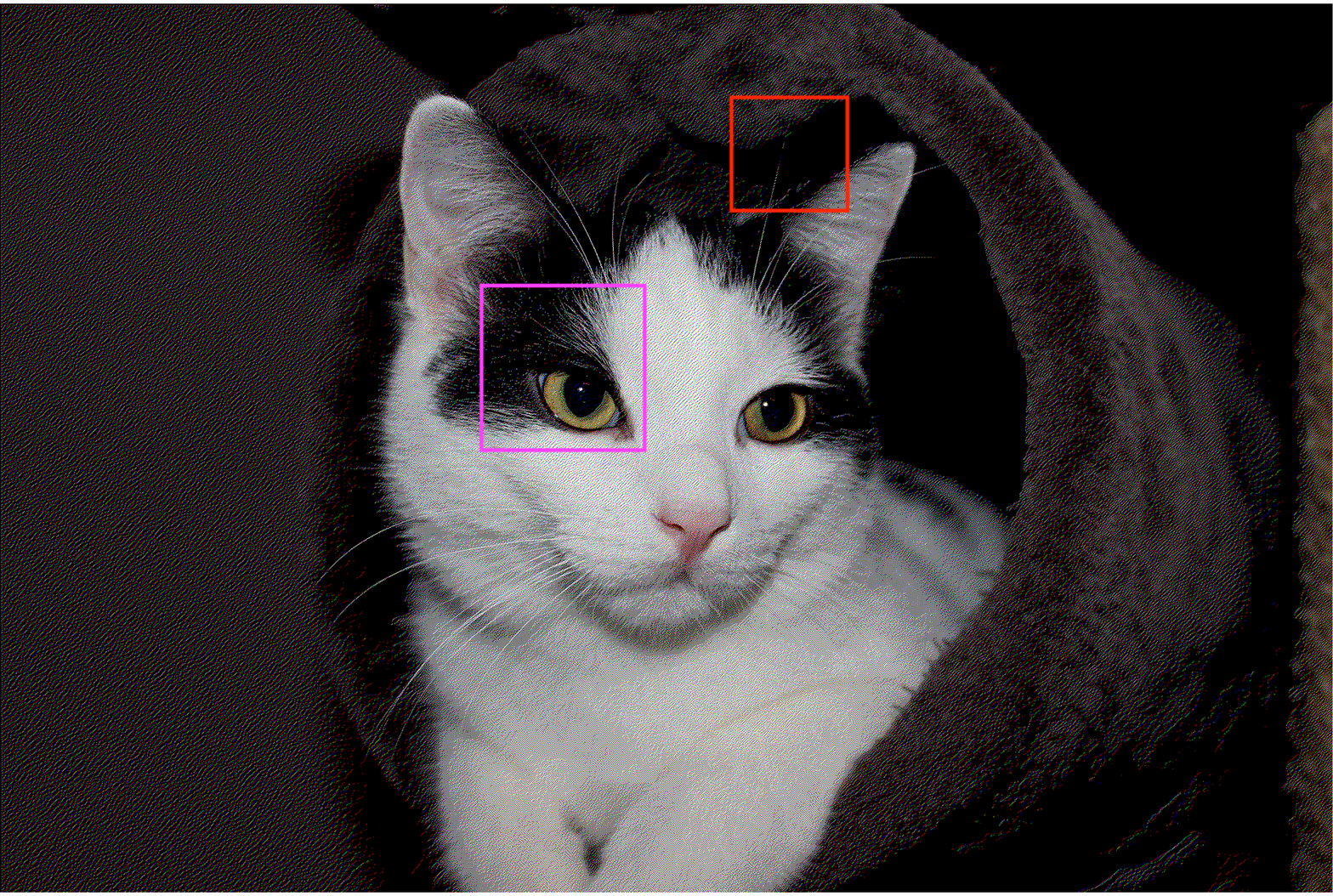}}%
}
\usebox{\bigpicturebox} \hfill
\begin{minipage}[b][\ht\bigpicturebox][s]{.20\textwidth}
\includegraphics[width=1.1\textwidth]{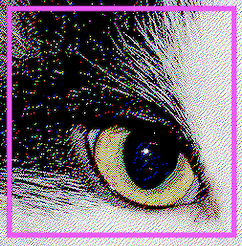}
\vfill
\includegraphics[width=1.1\textwidth]{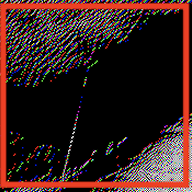}
\end{minipage}
\caption{Halftoned counterpart of image from Figure 6 produced using $1st$-$A$ scheme.
The value of the similarity index is $FSIM\!=\!0.9427$.}

\label{fig:10}
\end{figure}

\begin{figure}[ht!]
\centering
\sbox{\bigpicturebox}{%
  \scalebox{1}[1]{\includegraphics[width=.7\textwidth]{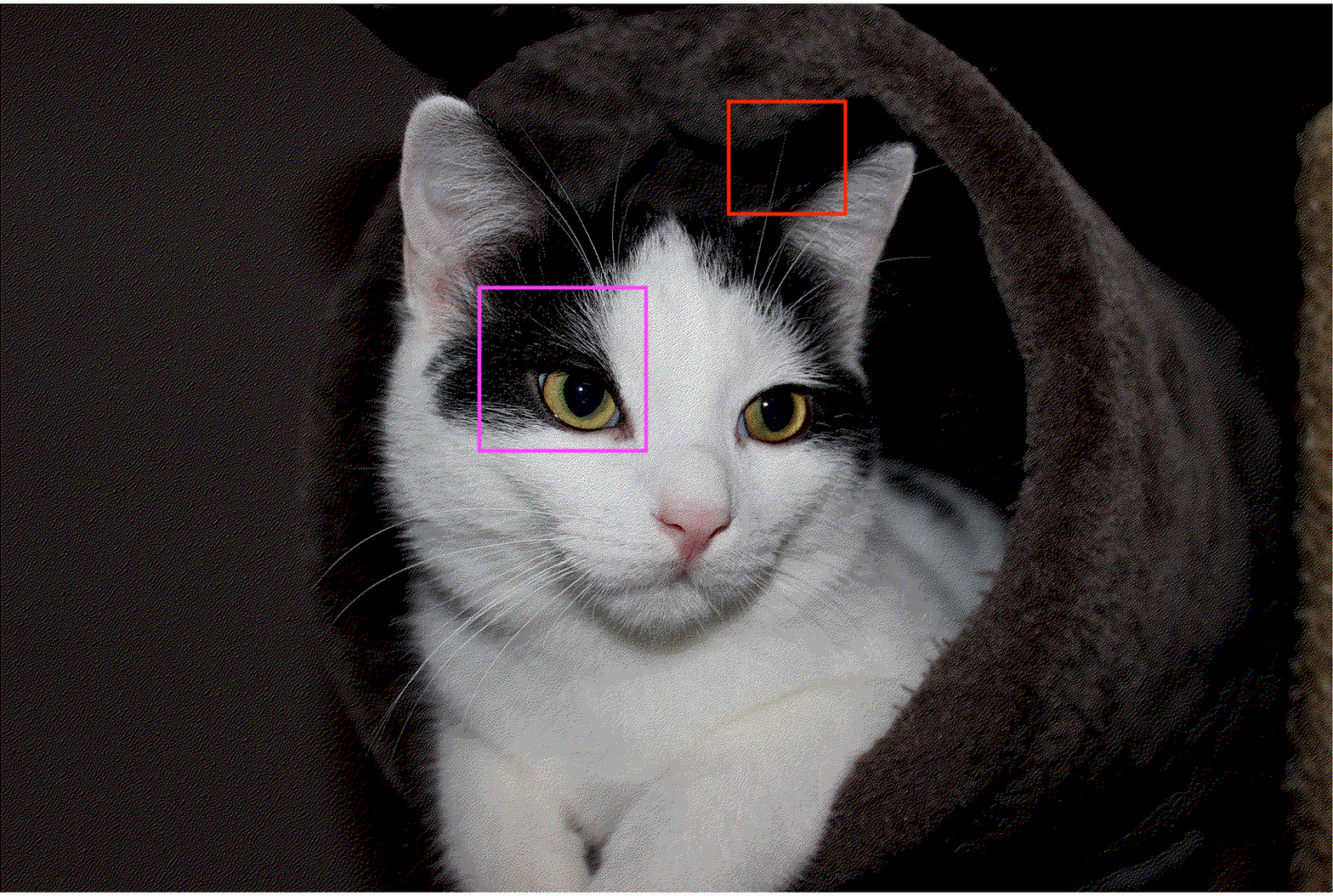}}}
\usebox{\bigpicturebox} 
\hfill
\begin{minipage}[b][\ht\bigpicturebox][s]{.20\textwidth}
\includegraphics[width=1.1\textwidth]{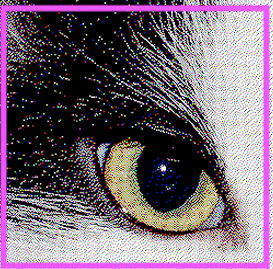}
\vfill
\includegraphics[width=1.1\textwidth]{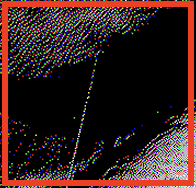}
\end{minipage}
\caption{Halftoned counterpart of image from Figure 6 produced using $2nd$-$SD$ scheme. The value of the similarity index is $FSIM\!=\!0.9511$.}

\label{fig:11}
\end{figure}
}

  Finally, we depict  one of the 50 images to illustrate  the visual quality of halftoned images arising from first and second order  weighted $\SD$ schemes.  Figure~\ref{fig:9} show the original color image and Figures~\ref{fig:10}, \ref{fig:11} depicts its halftoned counterparts. As one can  observe,  for the second order scheme, the halftoned image patterns are more refined and  the cat's whisker is completely reproduced, in contrast to the first order scheme. This is in line with our earlier observation that second order schemes  yield higher FSIM.

 \section{Discussion and Future Work}
  
In this paper, we proposed error diffusion algorithms  for digital halftoning based on $1$-bit weighted  $\SD$ quantization schemes. Even though these schemes are designed for best error decay in  the supremum norm, we observe excellent image quality also in terms of the commonly used {\it Feature Similarity Index}. 
Building on our findings, we see a number of interesting follow-up questions that we find worth investigating. 
First,  while the benefits of good reconstruction in the supremum norm seem to carry over to enhanced visual similarity to some extent, optimizing the weight matrices to minimize the former measure does not have the corresponding effect on the visual quality. This motivates  the question of whether  the weight matrices can be chosen directly to optimize structural similarity.  Second, in our numerical experiments, sparse filters of second order have been used, which have some performance limitations near the image boundary. In such situation, feedback filters chosen dynamically or based on the location may work better.  

 Lastly, the weighted $\SD$ quantization techniques developed in our current work  are especially designed for the two-dimensional image acquisition scenario.  Motivated by the growing importance of signal analysis on more sophisticated domains,  an important open question is how to generalize the concept of weighted $\SD$  schemes to higher dimensions as well as general manifold and graph geometries. 
 
\paragraph*{Acknowledgments.} The authors acknowledge support by the German Science Foundation (DFG) in the context of the Emmy Noether junior research group KR 4512/1-1 and the collaborative research center TR-109 as well as by the Munich Data Science Institute. FK would like to thank Sinan G\"unt\"urk for pointing out the connection between digital halftoning and $\SD$ quantization; furthermore the authors would like to thank Rongrong Wang for inspiring discussions related to the topic of this paper.


\section*{\Large Appendices}
\appendix
\addcontentsline{toc}{chapter}{Appendices}

\section{\large Taylor expansion}

For $\alpha=(\alpha_1,\alpha_2)$, $ \alpha_1, \alpha_2 \in \N$, we call  $|\al|=\al_1+\al_2\cdots +\al_d$ the order of $\al$. We  consider the factorial $\al!:= \al_1!\al_2!$, and write monomials as  $\x^{\al}= x_1^{\al_1}x_2^{\al_2}$ for $\x=(x_1,x_2) \in \R^2$.  
The partial derivative with respect to $\alpha$ of a function $f\!:\!\R^2\to \R$ is denoted by 
\begin{equation*}
 \da \al f= \da {\al_1}_1 \da {\al_2}_2 f= \frac{\da {|\al|} f}{\da {} x_1^{\al_1} \da {} x_2^{\al_2} } 
\end{equation*}


  Denote  by $C^k(\R^2)$ the class of bivariate functions $f$ for which all derivatives $\da \al f$ of order $|\al|\le k$ are continuous. Then we have the following Taylor approximation result.

\begin{thm}\cite{Hubbard2014}
Let $\mathcal{D}\subset\R^2$ be open, $f:\mathcal{D}\to \R$  be in $C^{k+1}(\R^2)$. Then for  $\bm a,\bm h \in \mathcal{D}$, such that  $\bm a+ \bm h\in \mathcal{D}$, $f$  can be represented as 
\begin{equation*}
    f(\bm a + \bm h)= \sum_{j=0}^k \sum\limits_{|\al|=j} \da {\al} f(\bm a) \frac{ \bb h ^{\al}}{\al!}+ R_{\bm a,k}(\bm h),
\end{equation*}
where the remainder $ R_{\bm a,k}$ is given in the integral form 

\begin{equation}\label{Taylor-rem-int-form}
    R_{\bm a,k}(\bm h):= (k+1)\sum\limits\frac{\bm h^\alpha}{\alpha!}\int\limits_{0}^1(1-t)^k\da {\al} f(\bm a+ t\bm h) \dx t.
\end{equation}
\end{thm}


\section{\large Estimation of Error Terms With Taylor Remainder}\label{Taylor-Rem-Est-1st-order}

As we have seen, a Taylor expansion of the quantization error for one-dimensional $\SD$ with a feedback filter $h$ satisfying  the moment conditions \eqref{moment-cond} and a  kernel $\Phi\in C^{r+1}(\R)$ gives rise to a combined remainder term of 
\begin{equation*}
 R_1:=\frac{1}{\lam} \sum_{n \in \N} v_n \sum\limits_{j=1}^L h_j R_{a_n, r}\Big(-\frac{j}{\lam}\Big)
\end{equation*}
where $a_n=x-\frac{n}{\lam}$  and the individual remainder terms $R_{a_n, r}$ are of the integral form 
\begin{equation*}
    R_{a_n,r}\big(-\!\frac{j}{\lam}\big)=\Big(\frac{-j}{\lam}\Big)^{r+1} \frac{1}{r!}\int\limits_0^1
(1-t)^r\Phi^{(r+1)}\big(a_n-\frac{t j}{\lam}\big) \dx t.
\end{equation*}
Then, an estimation of absolute value of $R_1$ leads to  
\begin{align*}
    |R_1|&=\frac{1}{\lam^{r+1} r!}\Big|\sum\limits_{j=1}^L h_j \cdot j^{r+1}\int_0^1 (1-t)^r \sum\limits_{n\in \N} \Phi^{(r+1)}\big(x-\frac{n}{\lam}-\frac{t j}{\lam}\big) \cdot\frac{v_n}{\lam}\,\dx t \Big|\\
    &\le \frac{\nofty{v}}{\lam^{r+1}}\sum\limits_{j=1}^L |h_j \cdot j^{r+1}| \int_0^1 |1-t|^r  \sum\limits_{n\in \N} \left|\Phi^{(r+1)}\big(x-\frac{n}{\lam}-\frac{t j}{\lam}\big)\right| \cdot\,\frac{1}{\lam}\, \dx t. 
\end{align*}
Denoting the constant 
$
\widetilde{C}_h:= \sum\limits_{j=1}^L |h_j| \cdot j^{r+1}
$
and using that $\int_0^1(1-t)^r\dx t= \frac{1}{r+1}$,
the range of $R_1$ can be upper-bounded as
\begin{equation*}
    |R_1|\le \frac{\nofty{v}}{\lam^{r+1}}\cdot \frac{\widetilde{C}_h}{(r+1)!}\cdot C\cdot \noone{\Phi^{(r+1)}}
\end{equation*}
for each value $x$, which shows that $|R_1|=\bO(\lam^{-(r+1)})$.




Analogously, for bivariate weighted $\SD$ schemes built from filters $h^{i,j}$ of the form discussed in  Section \ref{higher-order-SD-section} and a 
kernel $\Phi\in C^{r+1}(\R^n)$, the quantization error gives rise to a bivariate Taylor remainder which can be represented as 
\begin{align*}
 R_2 &:=  \frac{1}{\lam^2}  \sum\limits_{\bm n \in \N^2} v_{\bm n} \sum\limits_{i \, j} w_{i,j} \sum_{s=1}^L h^{i,j}_s R_{\bm a_{\bm n}, r} \big({-\tfrac{s}{\lam}\bb d_{i,j}}\big)  \\
    &=\sum\limits_{i \, j} w_{i,j} \sum_{s=1}^L h^{i,j}_s \tfrac{(-s)^{r+1}}{(r+1)^{-1}} \sum\limits_{|\alpha|=r+1}\tfrac{\bm d_{i j}^{\alpha}}{\alpha!\, \lam^{r+1}}\times \nonumber\\
   &\quad \quad \quad \quad \quad \quad \quad \quad \quad\times \int_0^1(1\!-\!t)^{r+1} \sum\limits_{\bm n \in \N^2} \da{\al} \Phi(\bm a_{\bm n}\!-\!\tfrac{ts}{\lam}\bb d_{i,j})\frac{v_{\bm n}}{\lam^2} \dx t. 
\end{align*}
 Observing that the sum inside of the integral can be bounded by a Riemann sum of $|\da \al \Phi(\cdot-\frac{ts}{\lambda}\bm d_{i,j})|$ and that
$\int_0^1(1-t)^r\dx t= \frac{1}{r+1}$, we obtain
\begin{equation}\label{talyor-part-est}
   |R_2|\le  \nofty{v}  \cdot C\cdot M_{\Phi} \cdot \widetilde{C}_h \cdot \sum\limits_{i \, j} w_{i,j} \sum\limits_{|\alpha|=r+1}\tfrac{|\bm d_{i j}|^{\alpha}}{\alpha! \lam^{r+1}},
\end{equation}
where  $ {M_{\Phi}\!:=\!\max\limits_{|\alpha|=r+1} \da \al\noone{\Phi} }$,
  $ {  \widetilde{C}_{h}\!:=\!\!\max\limits_{i\, j} \sum_{s=1}^L |h^{i,j}_s| s^{r+1}} $, and $C>0$ is a constant capturing the Riemann sum approximation error.  
By the definition $\bm d_{i,j}=\left(i,j \right)$, thus the two sums in \eqref{talyor-part-est} are equal to 

\vspace{-1mm}

\begin{equation}
   \sum\limits_{i \, j} w_{i,j} \sum\limits_{|\alpha|=r+1}\frac{|\bm d_{i j}|^{\alpha}}{\alpha! \lam^{r+1}}=
   \frac{1}{\lam^{r+1}}\cdot \widetilde{C}_{\bf W}
\end{equation}
with  $\widetilde{C}_{\bf W}:=\sum\limits_{i \, j} w_{i,j}\sum\limits_{m=0}^{r+1}\frac{i^{r+1-m}j^m}{(r+1-m)!m!}$ and one obtains that

\vspace{-1mm}

\begin{equation}
    |R_2|\le \frac{1}{\lam^{r+1}} \cdot  \nofty{v} \cdot \widetilde{C}_{\bf W} \cdot C \cdot \widetilde{C}_h  \cdot M_{\Phi}. 
\end{equation}


\begin{thebibliography}{99}
\fontsize{8}{5}\selectfont

\vspace{-1.1mm}

 \bibitem{Petersen1962} Petersen, D. P.,  Middleton, D. (1962). \emph{Sampling and reconstruction of wave-number-limited functions in N-dimensional Euclidean spaces.} Information and control, 5(4), 279-323.
 
\vspace{-1.1mm}


 
  \bibitem{Inose1963} Inose, H.,  Yasuda, Y. (1963). \emph{A unity bit coding method by negative feedback.} Proceedings of the IEEE, 51(11), 1524-1535.
  
  
\vspace{-1.1mm}


  \bibitem{DeFreitas1974} DeFreitas, R. (1974).\emph{ The low-cost way to send digital data: deltasigma modulation}. Electronic Design, 22, 68-73.
 
\vspace{-1.1mm}


 
  \bibitem{SteinbergFloyd}  Floyd R.W., Steinberg L. \emph{ An adaptive algorithm for spatial grey scale.} Proceedings of the Society of Information Display 17, 75–77 (1976).
  
  
\vspace{-1.1mm}


  
   \bibitem{Jarvis1976} Jarvis, J. F., Judice, C. N., and Ninke, W. H. (1976). \emph{A survey of techniques for the display of continuous tone pictures on bilevel displays.} Computer graphics and image processing, 5(1), 13-40.
  
\vspace{-1.1mm}

 
  \bibitem{Aziz1976}  Aziz, P. M., Sorensen, H. V., \& Vn der Spiegel, J. (1996).  \emph{An overview of sigma-delta converters.} IEEE signal processing magazine, 13(1), 61-84.
  
\vspace{-1.1mm}


   \bibitem{Knuth1987}  Knuth, D. E. (1987). \emph{Digital halftones by dot diffusion.} ACM Transactions on Graphics (TOG), 6(4), 245-273.

\vspace{-1.1mm}

 
  \bibitem{Knox1992}  Knox, K. T. (1992, May). \emph{Error image in error diffusion.} In Image Processing Algorithms and Techniques III (Vol. 1657, pp. 268-279). International Society for Optics and Photonics.
 
 
\vspace{-1.1mm}

   \bibitem{Kite1997}  Kite, T. D., Evans, B. L., Bovik, A. C.,  Sculley, T. L. (1997, October). \emph{Digital halftoning as 2-D delta-sigma modulation.} In Proceedings of International Conference on Image Processing (Vol. 1, pp. 799-802). IEEE.
   
  
\vspace{-1.1mm}
   
   
     \bibitem{Shiau1996} Shiau, J. N., \& Fan, Z. (1996, March). \emph{Set of easily implementable coefficients in error diffusion with reduced worm artifacts. }In Color Imaging: Device-Independent Color, Color Hard Copy, and Graphic Arts (Vol. 2658, pp. 222-225). International Society for Optics and Photonics. 
   
  
\vspace{-1.1mm}


   \bibitem{Ostromoukhov2001}  Ostromoukhov, V. \emph{A simple and efficient error-diffusion algorithm.} In E. Fiume, editor, Proc. SIGGRAPH 2001, Computer Graphics Proceedings, Annual Conference Series,
pages 567--572, Los Angeles, 2001.

\vspace{-1.1mm}


\bibitem{Calderbank2002}  Calderbank, A. R., \& Daubechies, I. (2002).\emph{ The pros and cons of democracy.} IEEE Transactions on Information Theory, 48(6), 1721-1725.
  
\vspace{-1.1mm}


    \bibitem{Secord2002} Secord, A. (2002, June). \emph{Weighted voronoi stippling.} In Proceedings of the 2nd international symposium on Non-photorealistic animation and rendering (pp. 37-43).
  
\vspace{-1.1mm}


 \bibitem{Kollig2003} Kollig, T., \& Keller, A. (2003). \emph{Efficient illumination by high dynamic range images}. In Proceed-
ings of the 14th Eurographics Workshop on Rendering, volume 44 of ACM International
Conference Proceeding Series, pages 45-50. 
  
  
\vspace{-1.1mm}


  
  \bibitem{Daubechies2003} Daubechies, I.,  DeVore, R. (2003). \emph{Approximating a bandlimited function using very coarsely quantized data: A family of stable sigma-delta quantizers of arbitrary order.} Annals of mathematics, 158(2), 679-710.
    
  
\vspace{-1.1mm}


  
  \bibitem{Gunturk2003} G\"unt\"urk, C. S. (2003). \emph{One‐bit sigma‐delta Quantization with exponential accuracy.} Communications on Pure and Applied Mathematics: A Journal Issued by the Courant Institute of Mathematical Sciences, 56(11), 1608-1630.
  

\vspace{-1.1mm}


\bibitem{Wang2004} Wang, Z., Bovik, A. C., Sheikh, H. R., \& Simoncelli, E. P. (2004).  \emph{Image quality assessment: from error visibility to structural similarity.} IEEE transactions on image processing, 13(4), 600-612.


\vspace{-1.1mm}


\bibitem{Schreier2005} Schreier R., \& Temes, G. C. (2005). \emph{Understanding delta-sigma data converters} (Vol. 74). Piscataway, NJ: IEEE press.

\vspace{-1.1mm}


\bibitem{Yilmaz2005}  Yilmaz, \"O. (2005). \emph{On coarse quantization of tight gabor frame expansions.} International Journal of Wavelets, Multiresolution and Information Processing, 3(02), 283-299.


\vspace{-1.1mm}

 \bibitem{Ostromoukhov2008} Vanderhaeghe, D., \& Ostromoukhov, V. (2008). \emph{Polyomino-based digital halftoning.} arXiv preprint arXiv:0812.1647.
 
\vspace{-1.1mm}


 
 \bibitem{Pang2008} Pang, W. M., Qu, Y., Wong, T. T., Cohen-Or, D., \& Heng, P. A. (2008). \emph{Structure-aware halftoning.} In ACM SIGGRAPH 2008 papers (pp. 1-8).
 
\vspace{-1.1mm}


 
 \bibitem{Krahmerthesis} Krahmer, F. (2009). \emph{Novel schemes for Sigma-Delta modulation: From improved exponential accuracy to low-complexity design (Doctoral dissertation, New York University)}. 

\vspace{-1.1mm}

 
 \bibitem{Balzer2009} Balzer, M., Schl\"omer, T., \& Deussen, O. (2009). \emph{Capacity-constrained point distributions: A variant of Lloyd's method.} ACM Transactions on Graphics (TOG), 28(3), 1-8.
 
 
\vspace{-1.1mm}



 \bibitem{Teuber2011}  Teuber, T., Steidl, G., Gwosdek, P., Schmaltz, C., \& Weickert, J. (2011). \emph{Dithering by differences of convex functions.} SIAM Journal on Imaging Sciences, 4(1), 79-108.
    
\vspace{-1.1mm}

   
\bibitem{Deift2011}   Deift, P., Krahmer, F., \& G\"unt\"urk, C. S. (2011).\emph{ An optimal family of exponentially accurate one‐bit Sigma‐Delta quantization schemes.} Communications on Pure and Applied Mathematics, 64(7), 883-919.
   
  \vspace{-1.1mm}
  
    \bibitem{Zhang2011} Zhang, L., Zhang, L., Mou, X., \& Zhang, D. (2011). \emph{ FSIM: A feature similarity index for image quality assessment.} IEEE transactions on Image Processing, 20(8), 2378-2386.
    
  
\vspace{-1.1mm}

  
       
     \bibitem{Krahmer2012} Krahmer, F., \& Ward, R. (2012). \emph{Lower bounds for the error decay incurred by coarse quantization schemes.} Applied and Computational Harmonic Analysis, 32(1), 131-138.  
 
\vspace{-1.1mm}

       
  \bibitem{Hubbard2014} Hubbard, J. H., \& Hubbard, B. B. (2015). \emph{Vector calculus, linear algebra, and differential forms: a unified approach} (pp. 818-pages). Matrix Editions.
  
\vspace{-1.1mm}

       
      \bibitem{Fornasier2016}   Fornasier, M., \& H\"utter, J. C. (2016). \emph{Consistency of probability measure quantization by means of power repulsion–attraction potentials.} Journal of Fourier Analysis and Applications, 22(3), 694-749.
  

\vspace{-1.1mm}

\bibitem{Wang2020} Lyu,\;H., Wang,\;R. (2020).\emph{ Sigma Delta quantization for images.} arXiv preprint arXiv:2005.08487.


 \end{thebibliography}
\end{document}